\documentclass[11pt,oneside]{article}
\usepackage{amssymb,amsmath}
\usepackage{url}
\usepackage[english]{babel}
\renewcommand{\le}{\leqslant}
\renewcommand{\ge}{\geqslant}
\newcommand{\bad}{\mathbf{Bad}}
\renewcommand{\L}{{\rm L}}
\newcommand{\area}{\mathbf{area}}

\newcommand{\RR}{\mathbb{R}}
\newcommand{\ZZ}{\mathbb{Z}}
\newcommand{\QQ}{\mathbb{Q}}

\newcommand{\NN}{\mathbb{N}}
\newcommand{\JJ}{\mathcal{J}}
\newcommand{\II}{\mathcal{I}}

\newcommand{\LL}{\mathbf{L}}
\newcommand{\CCC}{\mathcal{C}}
\newcommand{\KKK}{\mathcal{K}}
\newcommand{\RRR}{\mathcal{R}}
\newcommand{\SSS}{\mathcal{S}}
\newcommand{\vR}{\mathbf{R}}

\newcommand{\vr}{\mathbf{r}}
\newcommand{\vi}{\mathbf{i}}
\newcommand{\vf}{\mathbf{f}}
\newcommand{\maxab}{\max\{|A|^{1/i},|B|^{1/j}\}}
\newcommand{\maxhab}{\max\{|A|,|B|\}}

\newcommand{\f}{\mathbf{f}}
\newcommand{\bu}{{\bf u}}

\newtheorem{lemma}{Lemma}
\newtheorem{theorem}{Theorem}
\newtheorem{proposition}{Proposition}
\newtheorem{problem}{Conjecture}

\newtheorem{kg}{Theorem $\mathbf{1^{\prime}}$}

\newtheorem{kgg}{Theorem $\mathbf{2^{\prime}}$}

\newtheorem{corollary}{Corollary}

\newenvironment{proof}{\noindent \textsc{Proof.} \ }{\ \newline\hspace*{\fill}$\boxtimes$}

\newcommand{\cC}{{\cal C}}

\newcommand{\cM}{{\cal M}}

\newcommand{\cV}{{\cal V}}

\newcommand{\vA}{\mathbf{A}}

\newcommand{\vC}{C_0}

\begin{document}

\title{Badly approximable points on planar curves and a  \\
problem  of Davenport}

\author{
 Dzmitry Badziahin
\and
 Sanju Velani\footnote{Research partially supported by EPSRC grants EP/E061613/1 and EP/F027028/1  }
}

\date{{\small
{\em Dedicated  to our mathematical grandparents: \\ Harold
Davenport and Maurice Dodson  }}}

\maketitle

\begin{abstract}
Let  $\CCC$ be two times continuously differentiable  curve in
$\RR^2$  with  at least one point  at which the curvature is
non-zero. For any $i,j \ge 0$ with $i+j =1$, let  $\bad(i,j)$ denote
the set of  points $(x,y) \in \RR^2$ for which  $ \max \{
\|qx\|^{1/i}, \, \|qy\|^{1/j} \} > c/q $ for all $ q  \in \NN $.
Here $c = c(x,y)$ is a positive constant.  Our main result implies
that any finite intersection of such sets with $\CCC$ has full
Hausdorff dimension. This provides a solution to a problem  of
Davenport dating back to the sixties.
\end{abstract}

{\small


\medskip


}
\maketitle

\section{Introduction}

A real number $x$ is said to be {\em badly approximable} if there
exists a positive constant $c(x)$ such that
\begin{equation*}
 \|qx\| \ > \ c(x) \ q^{-1}  \quad \forall \  q \in \NN   \ .
\end{equation*}
Here and throughout $ \| \cdot  \| $ denotes the distance of a real
number to the nearest integer.  It is well known that set  $\bad $
of badly approximable numbers is of Lebesgue measure zero but of
maximal Hausdorff dimension; i.e. $ \dim \bad = 1 $. In higher
dimensions  there are various natural generalizations of $\bad$.
Restricting our attention to the plane $\RR^2$, given a pair of real
numbers $i$ and $j$ such that
\begin{equation}\label{neq1}
0\le i,j\le 1   \quad {\rm and \  } \quad  i+j=1  \, ,
\end{equation}
a point $(x,y) \in \RR^2$ is said to be {\em $(i,j)$-badly
approximable} if there exists a positive constant $c(x,y)$ such that
\begin{equation*}
  \max \{ \; \|qx\|^{1/i} \, , \ \|qy\|^{1/j} \,  \} \ > \
  c(x,y) \ q^{-1} \quad \forall \  q \in \NN   \ .
\end{equation*}
Denote by   $\bad(i,j)$  the set of $(i,j)$-badly approximable
points in $\RR^2$.  If $i=0$, then we use the convention that $
x^{1/i}\, :=0 $ and so $\bad(0,1)$ is identified with $\RR \times
\bad $. That is, $\bad(0,1)$ consists of points $(x,y)$ with $x \in
\RR$ and $y \in \bad$.  The roles of $x$ and $y$ are reversed if
$j=0$. In the case $i=j=1/2$, the set under consideration is the
standard set $\bad_2$  of simultaneously badly approximable points.
It easily follows from classical results in the theory of metric
Diophantine approximation   that $\bad(i,j)$ is of (two-dimensional)
Lebesgue measure zero and  it was shown in
\cite{Pollington-Velani-02:MR1911218} that $\dim \bad(i,j)=2$.

\subsection{The problem \label{prob}}
Badly approximable numbers obeying various functional relations were
first studied in the works of Cassels, Davenport and Schmidt from
the fifties and sixties.  In particular, Davenport
\cite{Davenport-64:MR0166154} in 1964 proved that for any $n \geq 2$
there is a continuum set of $\alpha\in\RR$ such that each of the numbers
$\alpha,\alpha^2,\dots,\alpha^n$ are all in $\bad$. In the same
paper,  Davenport \cite[p.52]{Davenport-64:MR0166154} states  ``{\em
Problems of a much more difficult character arise when the number of
independent parameters is less than the dimension of simultaneous
approximation. I do not know whether there is a set of $\alpha$ with
the cardinal of the continuum such that the pair $(\alpha,\alpha^2)$
is badly approximable for simultaneous approximation.}'' Thus, given
the  parabola $ \cV_2 := \{(x,x^2): x \in \RR\}$,  Davenport is
asking the question:
\begin{center}
{\em Is the set $\cV_2 \cap  \bad_2 $ uncountable?}
\end{center}
The goal of this paper is to answer this specific question for the
parabola and consider the  general setup involving an arbitrary
planar curve $\CCC $ and  $\bad(i,j)$.    Without loss of
generality, we assume that $\CCC $ is  given as a graph
$$
\CCC_f:= \{ (x,f(x))  : x \in I \}
$$
for some function $f$ defined on an interval $I \subset \RR$.  It is
easily seen that some restriction on the curve  is required to
ensure that $\CCC \cap  \bad(i,j)$ is not  empty.  For example,  let
$\L_{\alpha}$ denote the vertical line parallel to the $y$-axis
passing through the point $(\alpha,0)$ in the  $(x,y)$-plane.  Then,
it is  easily  verified, see
\cite[\S1.3]{Badziahin-Pollington-Velani-Schmidt} for the details,
that
$$
\L_\alpha \cap  \bad(i,j)  = \emptyset
$$
for any $\alpha \in \RR $ satisfying $ \liminf_{q \to \infty}
q^{1/i}  \|q \alpha\| = 0 \, . $ Note that the  $\liminf $ under
consideration  is  zero if $x$ is a Liouville number.  On the other
hand, if the $\liminf $ is strictly positive, which it is if $\alpha
\in \bad$, then
$$
\dim (\L_{\alpha} \cap  \bad(i,j))  = 1 \, .
$$
This result is much harder to prove and is at the heart of the proof
of Schmidt's Conjecture recently established in
\cite{Badziahin-Pollington-Velani-Schmidt}.    The upshot of this
discussion regarding vertical lines   is  that to build a general,
coherent theory for badly approximable points on planar curves we
need that the  curve $\CCC$ under consideration is in some sense
`genuinely curved'.  With this in mind, we will assume that  $\CCC$
is two times continuously differentiable and that there is at least
one point on $\CCC$  at which the curvature is non-zero.  We shall
refer to such a curve as a \emph{$C^{(2)}$ non-degenerate} planar
curve.  In other words and more formally, a planar curve $\CCC :=
\CCC_f$ is $C^{(2)}$ non-degenerate if $f \in C^{(2)}(I) $ and there
exits at least one point $ x \in I $ such that
$$ f''(x) \neq 0 \, .$$    For these curves,   it is reasonable to suspect that
$$
\dim (  \CCC  \cap \bad(i,j)) = 1 \ .
$$
%
%
%

\noindent If true,  this would imply that $\CCC \cap
\bad(i,j) $ is uncountable and since the parabola $\cV_2$  is a
$C^{(2)}$ non-degenerate  planar curve we obtain a positive answer
to Davenport's question.  To the best of our knowledge, there has
been no progress  with Davenport's question to date.  More
generally, for planar curves  (non-degenerate or not) the results
stated above for vertical lines constitute the first and essentially
only contribution.  The main result proved in this paper shows that
any finite intersection of $\bad(i,j)$ sets with a $C^{(2)}$
non-degenerate  planar curve  is of full dimension.

\subsection{The results \label{theresults} }

\begin{theorem}\label{thmfinite}
Let  $(i_1,j_1), \ldots ,(i_d,j_d)$ be a finite number of pairs of
real numbers  satisfying~\eqref{neq1}. Let $\CCC$ be a $C^{(2)}$
non-degenerate  planar curve. Then $$ \dim \Big(\bigcap_{t=1}^{d}
\bad(i_t,j_t)   \cap \CCC \Big) = 1 \ .
$$
\end{theorem}

A consequence of this theorem is the following statement regarding
the approximation of real numbers by algebraic numbers.  As usual,
the \emph{height} $H(\alpha)$ of an algebraic number is  the maximum
of the absolute values of the integer coefficients in its minimal
defining polynomial.

\begin{corollary}
 The set of $x\in\RR$ for which there exists a
positive  constant $c(x) $ such that
$$
|x-\alpha| >  c(x)   \, H(\alpha)^{-3}  \quad  \forall \ \mbox{real
algebraic numbers}  \ \alpha  \ \mbox{of degree} \, \leq 2
$$
is of full Hausdorff dimension.
\end{corollary}

\noindent The corollary  represents the  `quadratic' analogue of
Jarn\'{\i}k's  classical   $\dim \bad~=~1$  statement  and
complements the well approximable results of  Baker \& Schmidt
\cite{BakerSchmidt-1970}  and Davenport \& Schmidt
\cite{DavenportSchmidt-67}. It also makes a contribution to Problems
24, 25 and 26  in \cite[\S10.2]{Bugeaud-2004}.
To deduce the corollary from the theorem, we exploit the equivalent
dual form representation of the set $ \bad(i,j)$.  A  point
$(x,y)\in  \bad(i,j)$ if there exists a positive constant $c(x,y)$
such that
\begin{equation}\label{badij2}
\max\{|A|^{1/i}, |B|^{1/j}\}   \; \|Ax-By\|  >  c(x,y)  \qquad
\forall  \ (A,B)  \in\ZZ^2 \backslash \{(0,0)\}  \  .
\end{equation}
Then with $d=1$, $i=j=1/2$ and $\CCC = \cV_2$, the theorem implies that
$$
\dim \Big\{x\in\RR\;:\;  \max\{|A|^{2}, |B|^{2 }\}   \; \|Ax-Bx^2\|  >  c(x)  \
\forall  \ (A,B)  \in\ZZ^2 \backslash \{(0,0)\}\Big\}  = 1  \, .
$$
It can be verified that this is the  statement of the corollary
formulated in terms of integer polynomials.

\vspace{2ex}

Straight lines are an important class of $C^{(2)}$   planar curves
not covered by Theorem~\ref{thmfinite}. In view of the discussion in
\S\ref{prob},  this is to be expected since   the conclusion  of the
theorem is false for  lines in general.  Indeed,  it  is  only valid
for a  vertical line $\L_{\alpha} $  if  $\alpha$ satisfies the
Diophantine condition  $ \liminf_{q \to \infty}   q^{1/i}  \|q
\alpha\| > 0 \, . $  The following result provides an analogous
statement for  non-vertical lines.

\begin{theorem} \label{thmnvlines}
Let  $(i_1,j_1), \ldots ,(i_d,j_d)$ be a finite number of pairs of
real numbers  satisfying~\eqref{neq1}. Given $\alpha, \beta \in
\RR$, let $\L_{\alpha,\beta } $ denote the  line defined by the
equation $y = \alpha x +\beta $.   Suppose there exists $\epsilon>0$
such that
$$
\liminf_{q \to \infty}   q^{\frac{1}{\sigma}-\epsilon}  \|q \alpha
\|
> 0 \, \qquad{\rm where } \ \ \sigma := \max \{\min\{i_t,j_t\} : 1\leq t  \leq d
\} .
$$
Then
$$
\dim \Big(\bigcap_{t=1}^{d} \bad(i_t,j_t)   \cap \L_{\alpha,\beta } \Big) = 1 \ .
$$
\end{theorem}

\noindent In all likelihood this theorem is best possible apart from
the $\epsilon$ appearing in the Diophantine condition on the slope
$\alpha$ of the line. Indeed, this is the case for vertical lines --
see \cite[Theorem~2]{Badziahin-Pollington-Velani-Schmidt}.  Note that we always have that $\sigma \le 1/2$, so Theorem \ref{thmnvlines} is always valid for $\alpha \in \bad$. Also we point out that as a  consequence of the Jarn{\'\i}k-Besicovitch theorem, the  Hausdorff dimension of the exceptional set of $\alpha$ for which the conclusion of the theorem is not valid is bounded above by $2/3$.

\vspace{2ex}

%

\noindent{\em Remark 1. \ } The proofs of Theorem \ref{thmfinite}
and  Theorem \ref{thmnvlines} make use of a general Cantor framework
developed in \cite{Badziahin-Velani-MAD}.   The framework  is
essentially  extracted from the `raw' construction  used in
\cite{Badziahin-Pollington-Velani-Schmidt} to establish Schmidt's
Conjecture.  It will be apparent during the course of the proofs
that constructing the right type of general Cantor set in the  $d=1$
case is the main substance. Adapting the construction to
deal with finite intersections is not difficult and will  follow on
applying  the  explicit `finite intersection'  theorem  stated in
\cite{Badziahin-Velani-MAD}.  However, we  point out that  by
utilizing the arguments in
\cite[\S7.1]{Badziahin-Pollington-Velani-Schmidt} for countable
intersections   it is possible  to adapt the $d=1$ construction to
obtain the following strengthening of the theorems.

\begin{kg}\label{thmcountable}
Let  $(i_t,j_t)$ be a countable number of pairs of real numbers
satisfying~\eqref{neq1}  and  suppose that
\begin{equation} \label{liminfassump}
\liminf_{t \to \infty }  \min\{ i_t,j_t  \}  > 0   \ .
\end{equation}
Let $\CCC$ be a $C^{(2)}$ non-degenerate  planar curve.
Then
$$
\dim \Big(\bigcap_{t=1}^{\infty} \bad(i_t,j_t)   \cap \CCC
\Big) = 1 \ .
$$
\end{kg}

\begin{kgg}\label{thmcountable2}
Let  $(i_t,j_t)$ be a countable number of pairs of real numbers
satisfying~\eqref{neq1} and~\eqref{liminfassump}. Given $\alpha,
\beta \in \RR$, let $\L_{\alpha,\beta } $ denote the  line defined
by the equation $y = \alpha x +\beta $.   Suppose there exists
$\epsilon>0$ such that
$$
\liminf_{q \to \infty}   q^{\frac{1}{\sigma}-\epsilon}  \|q \alpha
\|
> 0 \, \qquad{\rm where } \ \  \sigma := \sup\{\min\{i_t,j_t\} : t
\in \NN \}.
$$
Then
$$
\dim \Big(\bigcap_{t=1}^{\infty} \bad(i_t,j_t)   \cap \L_{\alpha,\beta } \Big) = 1 \ .
$$
\end{kgg}

These statements should be true without the $\liminf$ condition
\eqref{liminfassump}. Indeed, without assuming \eqref{liminfassump}
the nifty argument developed by  Erez  Nesharim in \cite{Nesharim}
can be exploited to show  that the countable intersection of
the sets under consideration  are non-empty. Unfortunately, the
argument fails to show positive dimension let alone full dimension.

\vspace*{3ex}

\noindent{\em Remark 2.}  This
manuscript has taken a very   long time to produce. During its slow gestation, Jinpeng An \cite{Jinpeng}
circulated a paper  in which he shows that  $ \L_\alpha \cap
\bad(i,j)$ is winning (in the sense of Schmidt games -- see \cite[Chp.3]{Schmidt-1980}) for any
vertically line $ \L_{\alpha} $  with   $\alpha \in \RR$  satisfying
the Diophantine condition $ \liminf_{q \to \infty}   q^{1/i}  \|q
\alpha\| >  0 \, . $  An immediate consequence of this is that
$\bigcap_{t=1}^{\infty} \bad(i_t,j_t)   \cap \L_{\alpha} $  is of
full dimension as long as $\alpha $ satisfies the  Diophantine
condition  with $  i = \sup\{i_t : t \in \NN \} $.  The point is
that this is a statement  free of  \eqref{liminfassump} unlike  the
countable intersection result obtained in
\cite{Badziahin-Pollington-Velani-Schmidt}. In view of  An's  work   it is very tempting and
not at all outrageous to assert that {\em $ \bad(i,j)   \cap \CCC $
is winning} at least on the part   of the curve that is genuinely curved.   If true this would imply Theorem $1^{\prime}$ without
assuming \eqref{liminfassump}.  It is worth stressing that currently
we do not even know if $ \bad_{2}   \cap \CCC $ is winning.

\subsection{Davenport  in higher dimensions: what can we expect?}

For any  $n$-tuple of nonnegative real numbers $\vi:=(i_1, \ldots,i_n) $ satisfying $\sum_{s=1}^{n}  i_s = 1 $, denote by $\bad(\vi)$ the set of points $(x_1, \ldots,x_n) \in \RR^n$ for which there exists a positive constant $ c(x_1, \ldots,x_n)$ such that
\begin{equation*}
  \max \{ \; \|qx_1\|^{1/i_1} \,, \ldots  \ \|qx_n\|^{1/i_n} \,  \} \ > \
  c(x_1, \ldots,x_n) \ q^{-1} \quad \forall \  q \in \NN   \ .
\end{equation*}
The name of the game is to investigate the intersection of these
$n$-dimensional badly approximable sets with manifolds $ \cM \subset
\RR^n$.   A good starting point is to consider Davenport's problem
for arbitrary curves $\CCC$  in $\RR^n$. To this end and without
loss of generality, we assume that $\CCC $ is  given as a graph
$$
\CCC_{\vf} := \left\{ (f_1(x), \ldots, f_n(x))  : x \in I \right\}
$$
where  $\vf:= (f_1, \ldots, f_n) : I \to \RR^n $ is a map  defined
on an interval $I \subset \RR$. As in the planar case, to avoid
trivial empty intersection  with  $\bad(\vi)$ sets we assume that
the curve is genuinely curved. A curve $ \CCC :=\CCC_{\vf} \subset
\RR^n $ is said to be $C^{(n)}$  non-degenerate  if $ \vf \in
C^{(n)}(I) $ and there exists at least one  point $x \in I$  such
that the Wronskian
$$
w(f'_1, \ldots, f'_n ) (x) := \det(f_s^{(t)}(x))_{1\leq s,t \leq n}  \neq 0 \, .
$$
In the planar case ($n=2$),  this condition on the Wronskian is
precisely the same as saying that there exits at least one point on
the curve  at which the  curvature is non-zero. Armed with the
notion of  $C^{(n)}$  non-degenerate curves, there is no reason not
to believe in the truth of the following statements.

\begin{problem}
Let   $\vi_t :=(i_{1,t} \ldots,i_{n,t}) $ be a countable number of
$n$-tuples of non-negative real numbers satisfying $\sum_{s=1}^{n}
i_{s,t} = 1 $.  Let  $ \cC  \subset \RR^n $ be a $C^{(n)}$
non-degenerate curve. Then
$$
\dim \Big(\bigcap_{t=1}^{\infty} \bad(\vi_t)   \cap \cC
\Big) = 1 \ .
$$

\end{problem}

\begin{problem}
Let  $\vi:=(i_1, \ldots,i_n) $ be an $n$-tuple of non-negative real
numbers satisfying $\sum_{s=1}^{n}  i_s =~1 $. Let   $ \cC  \subset
\RR^n $ be a $C^{(n)}$  non-degenerate curve.  Then
 $
 \bad(\vi)   \cap \cC
$
is winning on some arc of $\cC$.
\end{problem}

\vspace*{3ex}

\noindent{\em Remark 1.}  In view of the fact that a winning set has
full dimension and that the intersection of countably many winning
sets is winning, it follows that Conjecture~B implies Conjecture~A.

\vspace*{2ex}

\noindent{\em Remark 2.}  Conjecture  A together with known
results/arguments from fractal geometry implies the strongest
version (arbitrary countable intersection plus full dimension) of
Schmidt's Conjecture in higher dimension:
$$\dim \Big(\bigcap_{t=1}^{\infty} \bad(\vi_t) \Big)  = n  \ . $$
In the case $n=2$,  this follows from An's result mentioned above
(Remark 2 in \S\ref{theresults})  -- see also his subsequent paper
\cite{Jinpeng2}.

\vspace*{2ex}

\noindent{\em Remark 3.}  Given that we basically know nothing in
dimension  $n > 2$,   a finite intersection version (including the
case $t=1$) of Conjecture  A would be a magnificent achievement.  In all
likelihood, any successful approach based on the general Cantor
framework developed in \cite{Badziahin-Velani-MAD} as in this paper
would yield Conjecture  A,  under the extra assumption involving the
natural analogue of the $\liminf$ condition \eqref{liminfassump}.

\vspace*{3ex}

We now turn our attention to general manifolds $ \cM   \subset
\RR^n$.  To avoid  trivial empty intersection  with
$\bad(\vi)$ sets, we assume that  the manifolds  under consideration
are non-degenerate. Essentially, these are smooth sub-manifolds of
$\RR^n$ which are sufficiently curved so as to deviate from any
hyperplane.  Formally, a  manifold ${\cal M}$ of dimension $m$
embedded in $\RR^n$ is said to be non-degenerate if it arises from a
non--degenerate map $\f:U\to \RR^n$ where $U$ is an open subset of
$\RR^m$ and ${\cal M}:=\f(U)$. The map $\f:U\to \RR^n:\bu\mapsto
\f(\bu)=(f_1(\bu),\dots,f_n(\bu))$ is said to be
\emph{non--degenerate at} $\bu\in U$ if there exists some $l\in\NN$
such that $\f$  is $l$ times continuously differentiable on some
sufficiently small ball centered at $\bu$ and the partial
derivatives of $\f$ at $\bu$ of orders up to $l$ span $\RR^n$.
\emph{If there exists at least one such non-degenerate point, we
shall say that the manifold $\cM=\f(U)$ is non--degenerate.}  Note
that in the case that the manifold is a curve $\cC$, this definition
is absolutely consistent with that of $\cC$  being $C^{(n)}$
non-degenerate. Also notice, that any
 real, connected analytic manifold not contained in any hyperplane of
$\RR^n$  is non--degenerate.  The following are the natural versions
of Conjectures A \& B for manifolds.

\begin{problem}
Let   $\vi_t :=(i_{1,t} \ldots,i_{n,t}) $ be a countable number of
$n$-tuples of non-negative real numbers satisfying $\sum_{s=1}^{n}
i_{s,t} =~1 $. Let  $ \cM  \subset \RR^n $ be a  non-degenerate
manifold. Then
$$
\dim \Big(\bigcap_{t=1}^{\infty} \bad(\vi_t)   \cap \cM
\Big) = \dim \cM \ .
$$
\end{problem}

\begin{problem}
Let  $\vi:=(i_1, \ldots,i_n) $ be an $n$-tuple of non-negative real
numbers satisfying $\sum_{s=1}^{n}  i_s =~1 $. Let  $ \cM  \subset
\RR^n $ be a  non-degenerate manifold.  Then
 $
 \bad(\vi)   \cap \cM
$
is winning on some patch of $\cM$.
\end{problem}


\vspace*{3ex}

\noindent{\em Remark 4.}  Conjecture  A together with the fibering
technique of Pyartly \cite{Pyartli-1969} should  establish Conjecture~C
for non-degenerate manifolds that can be foliated by non-degenerate
curves.  In particular, this  includes  any non-degenerate analytic
manifold
\footnote{A few days before completing this paper, Victor
Beresnevich communicated to us that he has established Conjecture~A under
the  extra assumption involving the natural analogue of
\eqref{liminfassump}.  In turn, under this  assumption,  by making
use of Pyartly's technique he  has proved  Conjecture~C for
non-degenerate analytic manifolds.  This in our opinion represents a
magnificent achievement  -- see Remark~3.}
.

\vspace*{3ex}

Beyond manifolds, it would be desirable  to investigate Davenport's
problem  within the more general context  of friendly measures
\cite{Kleinbock-Lindenstrauss-Weiss-04:MR2134453}.  We suspect  that
the above conjectures for manifolds remain valid with $\cM$ replaced by a
subset $X$ of $\RR^n$ that supports a  friendly measure.

\section{Preliminaries}

Concentrating on  Theorem \ref{thmfinite}, since any subset of a planar curve $\CCC$ is of dimension less than or equal to one we immediately obtain  that
\begin{equation}\label{klb} \dim \Big(\bigcap_{t=1}^{d}
\bad(i_t,j_t)   \cap \CCC \Big) \leq  1 \ .
\end{equation}
Thus, the proof of Theorem \ref{thmfinite}  reduces to establishing  the complementary lower bound statement  and as already mentioned  in \S\ref{thmfinite} (Remark 1) the crux is the $d=1$ case.  Without loss of generality, we assume that $i \le j$ .  Also, the case that $i=0$ is relatively straight forward  to handle  so let us assume that
\begin{equation}\label{neq1q}
0  <i \le j  <  1   \quad {\rm and \  } \quad  i+j=1  \, ,
\end{equation}
Then, formally the  key to establishing Theorem \ref{thmfinite} is the following statement.
\begin{theorem}\label{thmfinitelowerbd}
Let  $(i,j) $ be a  pair of
real numbers  satisfying~\eqref{neq1q}. Let $\CCC$ be a $C^{(2)}$
non-degenerate  planar curve. Then $$ \dim
\bad(i,j)   \cap \CCC  \geq 1 \ .
$$
\end{theorem}
The hypothesis that $\CCC=\CCC_f:= \{ (x,f(x))  : x \in I \}$ is  $C^{(2)}$
non-degenerate  implies that there exist
positive constants $\vC, c_0 > 0 $  so that
\begin{equation}\label{bound_curvat}
 c_0  \le |f'(x)| < \vC   \qquad {\rm and} \qquad c_0  \le |f''(x)| < \vC   \qquad  \forall  \ x \in I \, .
\end{equation}
To be precise, in general we can only guarantee  \eqref{bound_curvat} on a sufficiently small sub-interval $I_0$ of $I$.  Nevertheless, establishing Theorem \ref{thmfinitelowerbd} for the `shorter' curve  $\CCC^*_f = \{ (x,f(x))  : x \in I_0 \} $  corresponding to $f$ restricted to $I_0 $ clearly implies the desired dimension result for the curve~$\CCC_f $.

\medskip

To simplify notation the Vinogradov symbols $\ll$
and $\gg$ will  be used to indicate an inequality with an
unspecified positive multiplicative constant.  Unless stated otherwise,  the unspecified constant will at most be  dependant on $i,j,\vC$ and $c_0$ only.   If $a \ll b$ and $a
\gg b$ we write $ a \asymp b $, and  say that the quantities $a$
and $b$ are comparable.

\subsection{Geometric interpretation of $\bad(i,j)\cap\CCC$ }\label{subsec2_1}

We will work with the dual form of $\bad(i,j)$ consisting of points $(x,y) \in \RR^2$ satisfying  \eqref{badij2}.  In particular, for any
constant $c >0$, let  $\bad_{c}(i,j)$  denote  the set of points
$(x,y)\in \RR^2$  such that
\begin{equation}\label{neq2}
\max\{|A|^{1/i}, |B|^{1/j}\}   \; \|Ax-By\|  >  c   \qquad \forall \
(A,B)  \in\ZZ^2 \backslash \{(0,0)\}  \  .
\end{equation}

\noindent It is easily seen that $  \bad_{c}(i,j) \subset \bad(i,j)$
and
 $$\bad (i,j)    \, =   \, \bigcup_{c > 0}  \bad_{c}(i,j)    \ . $$
Geometrically,  given  integers  $A,B,C$  with $(A,B)\neq (0,0)$
consider the line $L=L(A,B,C)$ defined by the equation
$$
Ax-By+C = 0   \ .
$$
The set $\bad_{c}(i,j)$ simply consists of points in the plane that
avoid  the  $$ \frac{c}{ \max\{|A|^{1/i}, |B|^{1/j}\}  \;  \sqrt{A^2
+ B^2} }  $$ thickening  of each line $L$  -- alternatively, points
in the plane that lie within any such neighbourhood are removed. A
consequence of~\eqref{bound_curvat}  is that this thickening
intersects $\CCC$  in at most two closed arcs. Either of these arcs
will be denoted by $\Delta(L)$. Let $\RRR_0$ be the collection of
arcs $\Delta(L)$ on $\CCC$ arising from lines $L=L(A,B,C)$  with
integer coefficients and $(A,B)\neq (0,0)$.

The upshot of the above analysis is that the set
$\bad_c(i,j)\cap\CCC$ can be described as the set of all points
on~$\CCC$ that survive after removing the  arcs
$\Delta(L)\in\RRR_0$. Formally,
$$
\bad_c(i,j)\cap\CCC = \{(x,f(x))\in \CCC \, : \,  (x,f(x))\not\in
\Delta(L) \ \ \forall \, \Delta(L)\in \RRR_0 \}.
$$

\noindent For reasons that will become apparent later, it will be convenient to remove all but finitely many arcs.
With this in mind, let  $\SSS$ be a finite sub-collection of $\RRR_0$  and  consider the set
$$
\bad_{c,\SSS}(i,j)\cap\CCC = \{(x,f(x))\in \CCC \, : \,  (x,f(x))\not\in
\Delta(L) \ \ \forall \, \Delta(L)\in \RRR_0\backslash \SSS \}.
$$
Clearly, since we are removing fewer arcs $\bad_{c,\SSS}(i,j) \supset \bad_c(i,j)$. On the other hand,
$$
S:= \{(x,f(x))\in \CCC\;:\;  Ax-Bf(x)+C = 0 \ {\rm for  \ some  \ }  L(A,B,C) \ {\rm with \ }  \Delta(L) \in
\SSS  \}  \, $$   is a
finite set of points and it is easily verified  that
$$\bad_{c,\SSS}(i,j)\cap \CCC   \subset (\bad_c(i,j)\cap \CCC)
\cup S  \, . $$
Since $\dim S = 0$ for any finite set $S$ of points, Theorem \ref{thmfinitelowerbd} will follow on
showing that
\begin{equation} \label{show}
\dim \bad_{c,\SSS}(i,j)\cap \CCC\to 1\quad\mbox{ as }\quad c\to 0  \ .
\end{equation}

\noindent  In \S\ref{ofcourse} we will specify  exactly the finite collection of arcs $\SSS$ that are not to be  removed and  put  $\RRR := \RRR_0\backslash \SSS$  for this  choice of  $\SSS$.

\vspace*{2ex}

\noindent{\em Remark 1. } Without loss of generality,  when
considering lines $L=L(A,B,C)$ we will~assume that
\begin{equation*}\label{heighttriv}
(A,B,C)=1 \, .
\end{equation*}
 Otherwise we can divide  the
coefficients of $L$  by their common divisor. Then the
resulting line $L'$ will satisfy the required conditions and
moreover $\Delta(L') \supseteq \Delta(L)$.  Therefore, removing the arc $\Delta(L')$ from $\CCC$ takes care of removing
$\Delta(L)$.

\vspace*{2ex}

\medskip

\subsubsection{Working with the projection of $\bad_{c,\SSS}(i,j)\cap \CCC $ \label{fat}}

Recall that $\CCC = \CCC_f:= \{ (x,f(x))  : x \in I \}$ where $I
\subset \RR$ is an interval.  Let  $\bad^f_{c,\SSS}(i,j) $ denote the set
of $x \in I $ such that $ (x,f(x)) \in \bad_{c,\SSS}(i,j)\cap\CCC $.  In
other words  $\bad^f_{c,\SSS}(i,j) $ is the orthogonal projection of
$\bad_{c,\SSS}(i,j) \cap \CCC $ onto the $x$-axis. Now notice that in
view of ~\eqref{bound_curvat} the function $f$ is Lipschitz; i.e.
for some $\lambda >1$
$$
|f(x)-f(x')|\le \lambda |x-x'|\quad \forall x,x'\in I.
$$
Thus, the sets $\bad^f_{c,\SSS}(i,j) $ and $\bad_{c,\SSS}(i,j) \cap \CCC $ are
related by a bi-Lipschitz map and  so
$$
\dim \bad_{c,\SSS}(i,j) \cap \CCC   =  \dim \bad^f_{c,\SSS}(i,j) \, .
$$  Hence establishing \eqref{show}  is equivalent to showing that
\begin{equation} \label{show1}
\dim \bad^f_{c,\SSS}(i,j) \to 1\quad\mbox{ as }\quad c\to 0  \ .
\end{equation}
Next observe
that $ \bad^f_{c,\SSS}(i,j)$ can  equivalently be written as   the set of $x\in I$
such that $x\not\in \Pi(\Delta(L))$ for all $\Delta(L)\in\RRR_0\backslash \SSS$ where
the interval $\Pi(\Delta(L))\subset I$ is the orthogonal projection
of the arc $\Delta(L)\subset \CCC$ onto the $x$-axis.  Throughout the paper, we use the fact that the sets under consideration can be viewed either in terms of arcs $\Delta(L)$ on the curve $\CCC$  or sub-intervals  $\Pi(\Delta(L))$ of $ I$.   In order to  minimize unnecessary and cumbersome  notation,  we will simply write  $\Delta(L)$ even in the case of intervals and always refer to $\Delta(L)$ as an   interval.   It will be clear from the context whether $\Delta(L)$ is an arc on a curve or a genuine interval on $\RR$. However, we stress  that by the length of $\Delta(L)$ we
will always mean the length of the interval $\Pi(\Delta(L))$. In other words,
$$
|\Delta(L)| := |\Pi(\Delta(L))|.
$$

\subsection{An estimate for the size of  $\Delta(L)$}

Given a line $L=L(A,B,C)$, consider the function
$$
F_L\;:\;I\to \RR  \  : \  x  \to F_L(x):=Ax-Bf(x)+C.
$$
To simplify notation, if there is no risk of ambiguity  we shall  simply write $F(x)$ for $F_L(x)$. Now given an interval $\Delta(L)=\Delta(L(A,B,C))$ let
$$
V_L(\Delta):= \min_{x\in \Delta(L)}\{|F'_L(x)|\} = \min_{x\in \Delta(L)}\{|A-Bf'(x)|\}.
$$
Since $\Delta(L)$ is closed and $F_L$ is continuous the minimum always
exists. If there is no  risk of ambiguity  we shall  simply write $V_L$ for  $ V_L(\Delta) $. In short, the quantity $V_L$  plays a crucial role in estimating the size of $\Delta(L)$.

\begin{lemma}\label{lem_delta}
There exists an absolute constant $K \geq 1$ dependent  only on $i,j,\vC$ and $c_0$ such that
\begin{equation} \label{yeslabel}
|\Delta(L)|\le K\min\left\{\frac{c}{\maxab\cdot
V_L},\left(\frac{c}{\maxab\cdot |B|}\right)^{1/2}\right\}.
\end{equation}
\end{lemma}

\noindent \proof  The statement is essentially a  consequence of  Pyartly's~Lemma \cite{Pyartli-1969}: {\it Let
$\delta, \mu>0$ and $I\subset \RR$ be some interval. Let $f(x)\in
C^n(I)$ be function such that $|f^{(n)}(x)|>\delta$ for all $x\in
I$. Then there exists a contant $c(n)$ such that
$$
|\{x\in I\;;\; |f(x)|<\mu\}| \, \leq  \,  c(n) \,  \left(\frac{\mu}{\delta}\right)^{1/n}.
$$}

\noindent Armed with this, the first estimate for $|\Delta(L)|$  follows from the fact
that
$$
|F'_L(x)|\ge \delta:=V_L \quad\mbox{and}\quad |F_L(x)|\le \mu:=\frac{c}{\maxab}
$$
for all $x\in \Delta(L)$. The second makes use of the fact that
$$
|F_L''(x)| = |B f''(x)|>c_0 |B|\quad\forall  \;  x  \, \in \,  \Delta(L).
\vspace*{-4ex}
$$
\endproof

\vspace*{6ex}

\noindent{\em Remark 1. } The second term inside the minimum on the r.h.s. of \eqref{yeslabel} is absolutely crucial.
It shows that the length of $\Delta(L)$ can not be arbitrary
large  even when the quantity $V_L$ is small or even equal to zero.   The second term is not guaranteed if the curve is degenerate.  However,  for the lines (degenerate curves) $\L_{\alpha,\beta } $   considered in Theorem \ref{thmnvlines}   the Diophantine condition on $\alpha$  guarantees that  $V_L$ is not too small and hence allows us to adapt the   proof of Theorem \ref{thmfinitelowerbd} to this degenerate situation.

\subsubsection{Type 1 and Type 2  intervals \label{ofcourse} }

Consider an  interval $\Delta(L) = \Delta(L(A,B,C))\in \RRR$. Then
Lemma~\ref{lem_delta} implies  that
$$
\Delta(L)   \subseteq \Delta_1^*(L)    \qquad {\rm and } \qquad \Delta(L)   \subseteq \Delta_2^*(L)
$$ where the
intervals $\Delta_1^*(L)$ and $  \Delta_2^*(L)$ have the  same center as
$\Delta(L)$ and length given
$$
|\Delta_1^*(L)|:= \frac{2K\cdot c}{\maxab\cdot V_L} \; ,
$$
$$
|\Delta_2^*(L)|:= 2K\left(\frac{c}{\maxab\cdot |B|}\right)^{1/2}.
$$

\noindent We say that the  interval  $\Delta_1^*(L)$ is  of {\bf Type 1} and
$\Delta_2^*(L)$ is  of {\bf Type 2}. \emph{For obvious reasons, we assume that $B\neq 0$ in the case of Type 2}.  For each type of interval we define its
{\em height} in the following way:
$$
H(\Delta_1^*)= H(A,B) := c^{-1/2}\cdot V_L\cdot \maxab;
$$
$$
H(\Delta_2^*) = H(A,B) := (\maxab\cdot |B|)^{1/2}.
$$
So if   $\Delta^*(L)$ denotes an interval of either type we have that
$$|\Delta^*(L)| = 2Kc^{1/2}\cdot (H(\Delta^*))^{-1}.$$

\vspace*{1ex}

\noindent{\em Remark 1. }   Notice that for each positive number
$H_0$ there are only finitely many intervals $\Delta_2^*(L)$ of Type 2 such that $H(\Delta^*_2)\le H_0$.

\vspace*{4ex}

Recall,  geometrically $\bad_{c,\SSS}(i,j)\cap \CCC$ (resp. its
projection $ \bad^f_{c,\SSS}(i,j)$) is the set  of points on $\CCC$ (resp.
$I$) that survive after removing the intervals $\Delta(L) \in  \RRR_0\backslash \SSS$.
We now consider the corresponding subsets obtained by removing the larger intervals $\Delta^*(L)$.
Given $\Delta(L) \in \RRR_0$, the  criteria for which type of interval $\Delta^*(L)$ represents is
as follows.  Let $R \ge 2 $ be a large integer and  $\lambda$ be a  constant satisfying
\begin{equation}\label{cond_lambda}
\lambda>\max\left\{4, \frac1i, \frac{1+i}{j}\right\}.
\end{equation}
Furthermore, assume that the constant $c > 0$ satisfies
\begin{equation}\label{ineq_c}
c <\min\left\{(8(C_0+1)R^{-1-ij/2-\lambda})^2, ((C_0+1)C_0R^2)^{-2}\right\}.
\end{equation}
Given $\Delta(L)$ consider the associated Type 1 interval
$\Delta^*_1(L)$. There exists a unique $d\in\ZZ$ such that
\begin{equation}\label{cond_d}
R^d\le H(\Delta^*_1)<R^{d+1}.
\end{equation}
Choose $l_0$ to be the largest
integer such that
\begin{equation}\label{bound_l}
\lambda l_0\le\max\{d,0\}.
\end{equation}
Then we choose $\Delta^*(L)$ to be the interval $\Delta^*_1(L)$ of
Type 1 if
$$
V_L>(C_0+1)R^{-\lambda(l_0+1)}\max\{|A|,|B|\}.
$$
Otherwise, we take $\Delta^*(L)$ to be the interval $\Delta^*_2(L)$
of Type 2. Formally
\begin{equation} \label{yetnotanother}
\Delta^*(L):=\left\{\begin{array}{rl} \Delta_1^*(L)&\mbox{ if }\quad
V_L>(C_0+1)R^{-\lambda(l_0+1)}\max\{|A|,|B|\}.\\[2ex]
\Delta^*_2(L)&\mbox{ otherwise.} \end{array}\right.
\end{equation}

\vspace*{2ex}

\noindent{\em Remark 2. } It is easily  verified that for either type of interval, we have that
$$H(\Delta^*) \, \ge  \,
1.
$$
For Type 2 intervals  $\Delta^*_2(L)$ this follows by definition. For Type 1 intervals
$\Delta^*_1(L)$ assume that $H(\Delta_1)<1$. It then follows that $d<0$ and $l_0=0$.
In turn this implies that
\begin{eqnarray*}
H(\Delta_1) & := & c^{-1/2}V_L\maxab  \\[2ex] & \ge & c^{-1/2}
(C_0+1)R^{-\lambda}\max\{|A|,|B|\}\maxab  \\[1ex]
&\stackrel{\eqref{ineq_c}}\ge &  \max\{|A|,|B|\}\maxab \ \ge  \ 1 \, .
\end{eqnarray*}
This contradicts our assumption and thus we must have  that  $H(\Delta_1) \ge 1$.

\vspace*{3ex}

We now specify the finite sub-collection $\SSS$ of intervals from $\RRR_0$ which are not to be removed.
Let $n_0=n_0(c,R)$ be the minimal   positive integer satisfying
\begin{equation}\label{def_n0}
c^{1/2}\cdot R^{n_0}\cdot C_0\ge 1.
\end{equation}
Then,  define $\SSS$ to be  the collection of intervals  $\Delta(L) \in \RRR_0$ so that $\Delta^*(L)$ is
of Type 2 and   $H(\Delta^*)<R^{3n_0}$. Clearly  $\SSS$ is a finite collection of intervals -- see  Remark 1 above.  For this particular collection $\SSS$ we put
$$ \RRR :=   \RRR_0\backslash \SSS \, . $$

\vspace*{2ex}

Armed with this criteria  for choosing  $\Delta^*(L)$
given $\Delta(L)$  and indeed the finite collection $\SSS$ we consider the set
\begin{equation} \label{notanother}
\bad_c^*(i,j)\cap \CCC:=\{(x,f(x))\in \CCC\;:\; (x,f(x))\cap
\Delta^*(L)=\emptyset \;\;\forall  \;  \Delta(L)\in \RRR\}  \, .
\end{equation}
Clearly,
$$ \bad_c^*(i,j)\cap \CCC \subset
\bad_{c,\SSS}(i,j)\cap \CCC  \, $$ and so Theorem \ref{thmfinitelowerbd}
will follow on showing \eqref{show} with $\bad_{c,\SSS}(i,j)\cap \CCC$
replaced by $\bad_c^*(i,j)\cap \CCC$.
Indeed, from this point onward we will work with set defined by
\eqref{notanother}.  In view of this and to simplify notation we
shall simply redefine  $\bad_c(i,j)\cap \CCC$ to be
$\bad_c^*(i,j)\cap \CCC$  and write $\Delta(L) $ for $\Delta^*(L)$.
Just to make it absolutely clear,  the intervals  $\Delta(L)
:=\Delta^*(L)$ are determined  via the criteria
\eqref{yetnotanother} and $ \RRR $ is the collection of such
intervals arising from lines $L=L(A,B,C)$  apart from those associated with $\SSS$.  Also, the set  $\bad^f_c(i,j) $ is from this
point onward the orthogonal projection  of the redefined set
$\bad_c(i,j)\cap \CCC := \bad_c^*(i,j)\cap \CCC$.  With  this in
mind, the key to establishing \eqref{show1},   which in turn implies
\eqref{show} and therefore Theorem   \ref{thmfinitelowerbd}, lies in
constructing a Cantor-type subset $K_c(i,j)$ of $\bad^f_c(i,j)$
such that
\begin{equation*} \label{show2}
\dim K_c(i,j) \to 1\quad\mbox{ as }\quad c\to 0  \ .
\end{equation*}

\section{Cantor Sets and Applications }

The proof of Theorem \ref{thmfinite} and indeed Theorem
\ref{thmnvlines} makes use of a general Cantor framework developed
in \cite{Badziahin-Velani-MAD}.   This is what we now describe.


\subsection{A general Cantor framework \label{secgcf}}

\noindent{\bf The parameters.}  Let ${\rm I}$ be a closed interval
in $\RR$. Let $$\vR:=(R_n)  \quad {\rm with }  \quad {n\in \ZZ_{\ge
0}}$$ be a sequence of natural numbers and $$\vr:=(r_{m,n})   \quad
{\rm with }  \quad  m,n\in \ZZ_{\ge 0} \ {\rm  \ and \  }  \ m\le n
$$ be a two parameter sequence of non-negative real numbers.

\vspace*{2ex}

\noindent{\bf The  construction.}  We start by subdividing the
interval ${\rm I}$ into $R_0$ closed intervals $I_1$  of equal
length and denote by $\II_1$ the collection of such intervals.
Thus,
$$
\#\II_1  =  R_0   \qquad {\rm and } \qquad |I_1| =  R_0^{-1}\, |{\rm
I}|  \ .
$$
Next,  we remove  at most  $r_{0,0}$ intervals $I_1$ from $\II_1$ .
Note that we do not specify which intervals should be removed but
just give an upper bound on the number of intervals to be removed.
Denote by  $\JJ_1$ the resulting collection. Thus,
\begin{equation}\label{iona1}
\#\JJ_1  \ge    \#\II_1   -  r_{0,0}  \, .
\end{equation}
For obvious reasons, intervals in $\JJ_1$ will be referred to as
(level one) survivors.   It will be convenient to define  $\JJ_0 :=
\{J_0\} $ with $ J_0 :={\rm I} $.

\vspace*{1ex}

\noindent In general, for $n \ge 0$, given  a collection $\JJ_n$
we construct a nested collection $\JJ_{n+1}$  of closed intervals
$J_{n+1}$  using the following two operations.
\begin{itemize}
\item{\em Splitting procedure.} We subdivide each  interval $J_n\in \JJ_n$  into $R_n$ closed sub-intervals $I_{n+1}$ of equal length and denote by  $\II_{n+1}$ the collection of such intervals. Thus,
    $$
    \#\II_{n+1}  =  R_n \times  \#\JJ_n   \qquad {\rm and } \qquad |I_{n+1}| = R_n^{-1}\, |J_n|    \ .
    $$
\item{\em Removing procedure.} For each interval $J_n\in \JJ_n$ we remove at most
$r_{n,n}$ intervals $I_{n+1} \in  \II_{n+1} $ that lie  within $
J_n$.  Note that  the number of intervals $I_{n+1}$ removed is
allowed to  vary amongst  the  intervals  in  $\JJ_n$.    Let
$\II_{n+1}^{n} \subseteq \II_{n+1}  $ be the collection of
intervals that remain.  Next, for each interval $J_{n-1}\in
\JJ_{n-1}$ we remove at most  $r_{n-1,n}$  intervals $I_{n+1} \in
\II_{n+1}^{n} $ that lie within $ J_{n-1}$. Let  $\II_{n+1}^{n-1}
\subseteq \II_{n+1}^{n}  $ be the collection of  intervals that
remain. In general, for each interval $J_{n-k}\in \JJ_{n-k}$  $(1
\le k \le n)$ we remove  at most $r_{n-k,n}$ intervals $I_{n+1} \in
\II_{n+1}^{n-k+1} $ that lie within $J_{n-k}$. Also we let
$\II_{n+1}^{n-k} \subseteq \II_{n+1}^{n-k+1}  $ be the collection of
intervals that remain. In particular, $\JJ_{n+1}  := \II_{n+1}^{0} $
is the desired collection of (level $n+1$) survivors.  Thus, the
total number of intervals $I_{n+1}$ removed during the removal
procedure  is at most $
r_{n,n}\#\JJ_n+r_{n-1,n}\#\JJ_{n-1}+\ldots+r_{0,n}\#\JJ_0 $ and so
\begin{equation}\label{iona2}
\#\JJ_{n+1}\ge R_n\#\JJ_n-\sum_{k=0}^nr_{k,n}\#\JJ_k.
\end{equation}
\end{itemize}
\noindent Finally, having constructed the nested collections $\JJ_n$
of closed intervals   we consider the limit set
$$
 \KKK ({\rm I},\vR,\vr) :=  \bigcap_{n=1}^\infty \bigcup_{J\in
\JJ_n} J.
$$
The set $\KKK({\rm I},\vR,\vr)$ will be referred to as a {\em $({\rm
I},\vR,\vr)$ Cantor set.} For further details and examples see
\cite[\S2.2]{Badziahin-Velani-MAD}.  The following result
(\cite[Theorem 4]{Badziahin-Velani-MAD} enables us to estimate the
Hausdorff dimension of $\KKK({\rm I},\vR,\vr)$. It is the key to
establishing Theorem~\ref{thmfinite}.

\begin{theorem}\label{th_cantor}  Given $\KKK ({\rm I},\vR,\vr) $, suppose that
 $R_n\ge 4$ for all $n\in\ZZ_{\ge 0}$ and  that
\begin{equation}\label{cond_th2}
\sum_{k=0}^n \left(r_{n-k,n}\prod_{i=1}^k
\left(\frac{4}{R_{n-i}}\right)\right)\le \frac{R_n}{4}.
\end{equation}
Then
$$
\dim \KKK ({\rm I},\vR,\vr)  \ge \liminf_{n\to \infty}(1-\log_{R_n}
2).
$$
\end{theorem}

\noindent Here we use the convention that the product term in \eqref{cond_th2} is one when $k=0$  and by definition  $ \log_{R_n}\!2 := \log2/ \log R_n$.

\medskip

The next result \cite[Theorem 5]{Badziahin-Velani-MAD} enables us to
show that the intersection of finitely many sets $\KKK({\rm
I},\vR,\vr_i) $ is yet another $({\rm I},\vR,\vr) $ Cantor set for
some appropriately chosen $\vr$. This will enable us to establish
Theorem \ref{thmfinite}.

\begin{theorem}\label{th_cantor2}  For each integer $1\le i \le k $, suppose we are given a set $\KKK
({\rm I},\vR,\vr_i) $. Then
$$
\bigcap_{i=1}^{k} \KKK ({\rm I},\vR,\vr_i)   
$$ %
is a $({\rm I},\vR,\vr) $ Cantor set where
$$
\vr:=(r_{m,n})  \quad\mbox{with } \quad r_{m,n} :=\sum_{i=1}^k
r^{(i)}_{m,n}  \, .
$$
\end{theorem}

\subsection{The applications}\label{subsec3_2}

We wish to  construct an appropriate  Cantor-type set $K_c(i,j)\subset
\bad^f_c(i,j)$  which fits within  the general Cantor framework of \S\ref{secgcf}.   With this in mind, let $R \ge 2 $ be a large integer  and $$c_1:=c^{\frac12} R^{1+\omega}   \quad {\rm  where \ } \quad \omega:=\frac{ij}{4} $$
and the constant $c > 0$ satisfies \eqref{ineq_c}.
Take an interval $J_0\subset
I$ of length $c_1$.  With reference to \S\ref{secgcf} we
denote by $\JJ_0:=\{J_0\}$. We establish, by
induction on $n$, the existence of the collection $\JJ_n$ of closed
intervals $J_n$ such that $\JJ_n$ is nested in $\JJ_{n-1}$; that is, each interval $J_n$ in  $\JJ_n$  is contained in some interval  $J_{n-1}$ in  $\JJ_{n-1}$. The
length of an interval $J_n$ will be  given by
$$ |J_n |  :=  c_1R^{-n}   \,  , $$
and each
interval $J_n$ will satisfy the condition
\begin{equation}\label{cond_jn}
J_{n} \cap \Delta(L)
=\emptyset\qquad   \forall \;  L  \  \ {with}  \  \  H(\Delta)<R^{n-1} .
\end{equation}
In particular we put
$$
K_c(i,j):= \bigcap_{n=1}^\infty\bigcup_{J\in \JJ_{n}} J
$$
By construction, we have that
$$
K_c(i,j) \subset
\bad^f_c(i,j)
$$

%

\vspace*{2ex}

\noindent Now let $$\epsilon:= \frac{ijw}{2}=\frac{(ij)^2}{8} \qquad
{and}  \qquad  R>R_0(\epsilon) $$ be sufficiently large. Recall that
we are assuming that  $i>0$  and so $\epsilon $  is strictly
positive   -- we deal with  the   $i=0$  case  later in
\S\ref{pain}. Let $n_0= n_0(c,R)$ be the minimal  positive
integer satisfying \eqref{def_n0}; i.e.
$$
c^{1/2}\cdot R^{n_0}\cdot C_0\ge 1.
$$
It will be apparent from the construction of the collections of
$\JJ_n$ described in \S\ref{tyu}  that $K_c(i,j)$ is in fact a
$(J_0,\vR,\vr)$ Cantor set $ \KKK(J_0, \vR,\vr) $ with
$$\vR:=(R_n)=(R,R,R,\ldots)$$
and
$$\vr:=(r_{m,n})=\left\{\begin{array}{l}4R^{1-\epsilon}\quad\mbox{if }m=n;\\[2ex]
2R^{1-\epsilon}\quad\mbox{if }m<n,\; n-m\neq n_0\\[2ex]
3R^{1-\epsilon}\quad\mbox{if }n-m= n_0,\; n\ge 3n_0
\end{array}\right.
$$
By definition, note that for $R>R_0(\epsilon)$ large enough we have
that
$$
\mbox{l.h.s. of~\eqref{cond_th2} }=
\sum_{k=0}^nr_{n-k,n}\left(\frac{4}{R}\right)^k\le
4R^{1-\epsilon}\frac{1}{1-4/R}\le \frac{R}{4}=\mbox{r.h.s.
of~\eqref{cond_th2} }.
$$
Also note that  $R_n\ge 4$ for $R$ large enough. Then it  follows via
Theorem~\ref{th_cantor} that
$$
\dim \bad^f_c(i,j) \ge \dim K_c(i,j) =   \dim \KKK(J_0, \vR,\vr) \ge 1-\log_R 2 \, .
$$
This is true for all $R$ large enough (equivalently all $c>0$ small enough) and so   on letting  $R\to\infty$ we obtain  that
$$
\dim \bad(i,j)\cap\CCC  \geq  \dim \bad^f_c(i,j)  \to
1.
$$
This  proves Theorem \ref{thmfinitelowerbd} modulo the construction of the collections $\JJ_n$.  Moreover,  Theorem~\ref{th_cantor2} implies that
$$
\bigcap_{t=1}^d (\bad(i_t,j_t)) \cap \CCC
$$
contains the Cantor-type set $\KKK(J_0,\vR,\tilde{\vr})$  with
$$\tilde{\vr}:=(\tilde{r}_{m,n})=\left\{\begin{array}{l}4dR^{1-\tilde{\epsilon}}\quad\mbox{if }m=n;\\[2ex]
2dR^{1-\tilde{\epsilon}}\quad\mbox{if }m<n,\; n-m\neq n_0\\[2ex]
3dR^{1-\tilde{\epsilon}}\quad\mbox{if }n-m=n_0,\; n\ge 3n_0.
\end{array}\right.
$$
where
$$
\tilde{\epsilon} := \min_{1\le t\le d}\left(\frac{(i_tj_t)^2}{8}\right).
$$
On applying Theorem~\ref{th_cantor} to the set
$\KKK(J_0,\vR,\tilde{\vr})$ and letting $R\to\infty$ implies that
$$\dim \Big(\bigcap_{t=1}^{d}
\bad(i_t,j_t)   \cap \CCC \Big) \ge  1 \, .  $$
This together with the upper bound statement  \eqref{klb} establishes
Theorem~\ref{thmfinite} modulo of course  the construction of the collections $\JJ_n$ and the assumption that $i > 0$.

\section{Preliminaries for constructing $\JJ_n$ \label{sec_split}}

In order to construct the  appropriate collections $\JJ_n$  
described in \S\ref{subsec3_2}, it is necessary to partition
the collection  $\RRR$ of intervals $\Delta(L) $  into various classes.  The aim  is to have  sufficiently good control on the parameters $|A|, |B|$ and $V_L$ within each class.  Throughout, $R \ge 2 $ is  a large integer.

\vspace*{3ex}

\noindent $\bullet$ \emph{Firstly we partition  all Type 1 intervals $\Delta(L)\in\RRR$  into
classes $C(n)$ and $C(n,k,l)$.}

\vspace*{0.3ex}

A Type 1 interval  $\Delta(L)\in C(n)$ if
\begin{equation}\label{classmore}
R^{n-1} \, \le \,  H(\Delta)  \, <  \, R^{n}    \, .
\end{equation}

\noindent Furthermore, $\Delta(L)\in C(n,k,l) \subset C(n)  $ if
\begin{equation}\label{class_prop1}
2^kR^{n-1}\le H(\Delta)<2^{k+1}R^{n-1} \qquad 0\le k< \log_2 \! R
\end{equation}
and
\begin{equation}\label{class_prop2}
R^{-\lambda (l+1)}(\vC+1)\max\{|A|,|B|\}<V_L\le R^{-\lambda
l}(\vC+1)\max\{|A|,|B|\}  \, .
\end{equation}

\noindent Note that since the intervals $\Delta(L)$ are of Type 1, it follows from~\eqref{bound_l} that
$l \le l_0$.  Moreover
$$
V_L=|A-Bf'(x_0)|\stackrel{\eqref{bound_curvat}}\le |A|+C_0|B|\le
(1+C_0)\max\{|A|,|B|\}
$$
so $l$ is also nonnegative.  Here and throughout  $x_0$ is the point at which $|F_L'(x)| = |A-Bf'(x)| $ attains its minimum with $x\in \Delta(L)$.  We let
$$
C(n,l):=\bigcup_{k=0}^{\log_2 \! R} C(n,k,l).
$$

\noindent $\bullet$  \emph{Secondly we partition  all Type 2 intervals $\Delta(L)\in\RRR$  into
classes $C^*(n)$ and $C^*(n,k)$.}

\vspace*{0.3ex}

 A Type 2 interval $\Delta(L)\in C^*(n)$ if~\eqref{classmore} is satisfied. Furthermore, $\Delta(L)\in C^*(n,k) \subset C^*(n)$ if~\eqref{class_prop1} is satisfied.

\medskip

\noindent Note that since   $H(\Delta)\ge 1$, we have the following the complete split of $\RRR$:
$$
\RRR = \left(\bigcup_{n=0}^\infty C(n)\right) \cup\left(\bigcup_{n=0}^\infty
C^*(n)\right).
$$
We now investigate the  consequences of the above classes on the parameters $|A|, |B|$ and $V_L$ and introduce further subclasses to gain tighter control.

\subsection{Estimates for $|A|$, $|B|$ and $V_L$ within
a given class}

\subsubsection{Class $C(n,k,l)$ with $l\ge 1$ \label{dbq}}

Suppose $\Delta(L(A,B,C))\in C(n,k,l)$ for some  $l\ge 1$. By definition each
of these classes corresponds to the case that the derivative $V_L =
|F'_L(x_0)|$ satisfies~\eqref{class_prop2}. In other words the derivative is essentially smaller than the expected value $\max\{|A|,|B|\}$.  Now observe that the  r.h.s. of \eqref{class_prop2} implies either
$$
|A-f'(x_0) B|<\frac{\vC+1}{R^\lambda}|A| \Leftrightarrow
\left(1-\frac{\vC+1}{R^\lambda}\right)< \frac{|f'(x_0)
B|}{|A|} <\left(1+\frac{\vC+1}{R^\lambda}\right)
$$
or
$$
|A-f'(x_0) B|<\frac{\vC+1}{R^\lambda}|B|   \Leftrightarrow\!
\left(1-\frac{\vC+1}{|f'(x_0)|R^\lambda}\right) <\frac{|A|}{|f'(x_0) B|} <\!\!\left(1+\frac{\vC+1}{|f'(x_0)|R^\lambda}\right)\, .
$$
Since $|f'(x_0)|\ge c_0> 0$ then in both cases,  for $R$ large enough
we have that
\begin{equation}\label{class_eq1}
2^{-1} |A|<|f'(x_0) B|< 2|A|\quad\mbox{or} \quad |A|\asymp |B|.
\end{equation}

\noindent On substituting the estimate~\eqref{class_prop2} for $V_L$ into the definition of the height $H(\Delta)$  we obtain  that
$$
c^{-\frac12}\cdot|A|^{\max\{\frac{i+1}{i},\frac{j+1}{j}\}}R^{-\lambda
(l+1)}\ll H(\Delta)\ll
c^{-\frac12}\cdot|A|^{\max\{\frac{i+1}{i},\frac{j+1}{j}\}}R^{-\lambda
l}.
$$
This together with \eqref{class_prop1} and the fact that $i\le j$, implies that
\begin{equation}\label{class_eq2}
\left(\frac{2^kc^{\frac12}}{R} R^{n+\lambda
l}\right)^{\frac{i}{i+1}}\ll |A|,|B|\ll
\left(\frac{2^{k}c^{\frac12}}{R}
R^{n+\lambda(l+1)}\right)^{\frac{i}{i+1}}.
\end{equation}

%

\subsubsection{Class $C(n,k,0)$ \label{time}}

By \eqref{class_prop1} and \eqref{class_prop2},  we have that in this
case
$$
c^{-\frac12}\cdot\frac{\max\{|A|,|B|\}}{R^\lambda}\max\{|A|^{1/i},|B|^{1/j}\}\ll
H(\Delta)\ll \frac{2^k}{R}R^n.
$$
Therefore,

\begin{equation}\label{class0_aeq}
|A|\ll
\left(\frac{2^kc^{\frac12}}{R}R^{n+\lambda}\right)^{\frac{i}{i+1}}
\end{equation}
\begin{equation}\label{class0_beq}
|B|\ll
\left(\frac{2^kc^{\frac12}}{R}R^{n+\lambda}\right)^{\frac{j}{j+1}}.
\end{equation}

\noindent Unfortunately these bounds for $|A|$ and $|B|$ are not strong  enough for our purpose.
Thus, we   partition  the class $C(n,k,0)$ into the following subclasses:
\begin{enumerate}
\item[] $C_1(n,k):=\{\Delta(L(A,B,C))\in C(n,k,0)\;:\; |A|\ge\frac12 |f'(x_0)||B|\}$

\item[] $C_2(n,k):=\{\Delta(L(A,B,C))\in C(n,k,0)\;:\; |A|<\frac12|f'(x_0)||B|,\; |A|^{1/i}\le |B|^{1/j}\}$

\item[] $C_3(n,k):=\{\Delta(L(A,B,C))\in C(n,k,0)\;:\; |A|<\frac12|f'(x_0)||B|,\; |A|^{1/i}>
|B|^{1/j}\}$.
\end{enumerate}

\medskip

\noindent \noindent $\bullet$  {\bf Subclass $C_1(n,k)$ of $C(n,k,0)$.} By \eqref{class0_aeq} we have
the following bounds for $|B|$ and $V_L$:
\begin{equation}\label{class1}
|B| \, , \  V_L\ll
\left(\frac{2^kc^{\frac12}}{R}R^{n+\lambda}\right)^{\frac{i}{i+1}}  \ .
\end{equation}
Note that this bound for $ |B|$ is  stronger than~\eqref{class0_beq}.

\medskip

\noindent  \noindent $\bullet$  {\bf Subclass $C_2(n,k)$ of $C(n,k,0)$.} We can strengthen the bound \eqref{class0_aeq} for
$|A|$ by the following:
\begin{equation}\label{class2_aeq}
|A|\le |B|^{i/j}\ll
\left(\frac{2^kc^{\frac12}}{R}R^{n+\lambda}\right)^{\frac{i}{j+1}}.
\end{equation}
Since $|A|<\frac12 |f'(x_0)||B|$ we have that $V_L\asymp|B|$,
therefore
\begin{equation}\label{class2_veq}
 V_L\ll
\left(\frac{2^kc^{\frac12}}{R}R^{n+\lambda}\right)^{\frac{j}{j+1}}.
\end{equation}

\noindent Also  we get that $\maxab= |B|^{1/j}$ which together
with~\eqref{class_prop1} implies that for any two
$\Delta(L_1(A_1,B_1,C_1)), \Delta(L_2(A_2,B_2,C_2))\in C_2(n,k)$,
\begin{equation}\label{class2_b12}
V_{L_1}\asymp B_1\asymp B_2\asymp V_{L_2}.
\end{equation}

\medskip

 \noindent $\bullet$  {\bf Subclass $C_3(n,k)$ of $C(n,k,0)$.} As with the previous subclass $C_2(n,k)$ we
have that
$$V_L\asymp |B|    \qquad \forall \quad   \Delta(L(A_2,B_2,C_2))\in C_3(n,k) \, . $$
We  partition  $C_3(n,k)$ into subclasses
$C_3(n,k,u,v)$ consisting of intervals $\Delta(L(A,B,C))\in C_3(n,k)$ with
\begin{equation}\label{class3_prop}
2^v R^{\lambda u}|B|^{1/j}< |A|^{1/i}\le 2^{v+1}R^{\lambda
u}|B|^{1/j} \qquad u\ge 0 \qquad \lambda\log_2 \! R \ge v\ge 0.
\end{equation}
Then
$$
|B|^{\frac{j+1}{j}}R^{\lambda u}< |B||A|^{1/i}\asymp
V_L\maxab = c^{\frac12} H(\Delta) \stackrel{\eqref{class_prop1}}<\frac{2^{k+1}c^{\frac12}}{R}R^n.
$$
Therefore
\begin{equation}\label{class3_beq}
V_L\asymp |B|\ll
\left(\frac{2^kc^{\frac12}}{R}R^{n}\right)^{\frac{j}{j+1}}R^{-\frac{\lambda
uj}{j+1}}
\end{equation}
and
\begin{equation}\label{class3_aeq}
|A|\stackrel{\eqref{class3_prop}}\ll R^{\lambda (u+1)i}|B|^{i/j}\ll
\left(\frac{2^kc^{\frac12}}{R}R^{n}\right)^{\frac{i}{j+1}}R^{\frac{\lambda
uij}{j+1}+\lambda i}.
\end{equation}

We proceed with estimating  the size of the parameter $u$. The fact that $|A|<\frac12
|f'(x_0)||B|$ together with~\eqref{class3_prop} and~\eqref{class3_beq} implies that
$$
R^{\lambda
u}\stackrel{\eqref{class3_prop}}<\frac{|A|^{1/i}}{|B|^{1/j}}\ll
|B|^{\frac{j-i}{ij}}\stackrel{\eqref{class3_beq}}\ll
R^{\frac{(j-i)n}{i(j+1)}}.
$$
Therefore for large $R$, if $C_3(n,k,u,v)$ is nonempty then
$u$ satisfies
\begin{equation}\label{class3_lbound}
0\le \lambda u\le \frac{j-i}{i(1+j)}\cdot n+1.
\end{equation}
In particular, this  shows that  $u$ is
smaller than $n$ if  $\lambda>1/i$.
Finally,  it can be verified  that the inequalities given by~\eqref{class2_b12} are valid for any two intervals $\Delta(L_1),
\Delta(L_2)\in C_3(n,k,u,v)$.

\subsubsection{Class $C^*(n,k)$}

By the definition~\eqref{bound_l} of $l_0$, we have that
$$V_L \le  R^{-\lambda}(\vC+1)\max\{|A|,|B|\}  \, . $$
This corresponds to the  r.h.s. of  \eqref{class_prop2} with $l=1$  and thus the same arguments as in  \S\ref{dbq} can be utilized to  show that ~\eqref{class_eq1} is satisfied. By substituting this into the definition of the height we obtain that
$$
H(\Delta)\asymp |A|^{\frac{i+1}{2i}}
$$
which in view of~\eqref{class_prop1} implies that
\begin{equation}\label{class_eq2_2}
|A|\asymp |B|\asymp \left(\frac{2^k}{R}\cdot
R^n\right)^{\frac{2i}{i+1}}.
\end{equation}

\noindent A consequence of this estimate is that all intervals
$\Delta(L)\in C^*(n,k)$ have comparable coefficients $A$ and $B$. In
other words, if $\Delta(L_1),\Delta(L_2)\in C^*(n,k)$ then
$$|A_1|\asymp|B_1|\asymp|B_2|\asymp|A_2| \, .$$
To estimate the size of $V_L$ we make use of the fact that
\begin{eqnarray}
V_L&\le &(C_0+1)R^{-\lambda(l_0+1)}\max\{|A|,|B|\}\stackrel{\eqref{bound_l}}\le
(C_0+1)R^{-d}\max\{|A|,|B|\}  \nonumber\\[2ex]
&\stackrel{\eqref{cond_d}}\le
&(C_0+1)c^{1/2}\cdot\frac{\max\{|A|,|B|\}}{V_L\maxab}\nonumber
\end{eqnarray}
This together with~\eqref{ineq_c} and \eqref{class_eq2_2} enables us to verify that
\begin{equation}\label{type2_eq}
V_L\le \frac{|B|}{R\cdot H(\Delta)}\ll \left(\frac{2^k}{R}\cdot
R^n\right)^{-\frac{j}{i+1}}.
\end{equation}

\subsection{Additional subclasses $C(n,k,l,m)$ of  $C(n,k,l)$ \label{additional}}

It is necessary  to partition each  class  $C(n,k,l)$ of Type 1 intervals $\Delta(L)$ into the following subclasses
to provide stronger control on  $V_L$.   For $m\in \ZZ$,  let
\begin{equation}\label{class0m_prop}
C(n,k,l,m)\!:=\!\left\{\!\!\Delta(L(A,B,C))\!\in C(n,k,l)\left|\!\!
\begin{array}{l}
2^{-m-1}R^{-\lambda l}(\vC+1)\max\{|A|,|B|\}<\!V_L\\[1ex]
V_L\le 2^{-m}R^{-\lambda l}(\vC+1)\max\{|A|,|B|\}
\end{array}
\right.\!\!\!\!\right\}.
\end{equation}
In view of \eqref{class_prop2}, it is easily verified that
$$ 0   \leq m \leq \lambda\log_2 \! R\asymp \log R  \, . $$
An important consequence  of  introducing these subclasses  is  that for any  two
intervals $\Delta(L_1),\Delta(L_2)$ from $C(n,k,l,m) $ with $l\ge 1$ or
from $C_1(n,k)\cap C(n,k,0,m)$, we have that
\begin{equation}\label{asymp_prop}
V_{L_1}\asymp V_{L_2}\quad\mbox{ and }\quad |A_1|\asymp |A_2|.
\end{equation}



\section{Defining the collection $\JJ_n$ \label{tyu}}

We  describe the procedure for constructing the
collections $\JJ_n \ (n = 0,1,2, \ldots)$  that lie at the heart of the construction of the Cantor-type set $ K_c(i,j) =\KKK(J_0, \vR,\vr)$ of \S\ref{subsec3_2}.  Recall that  each interval $J_n \in  \JJ_n$  is to be nested in some interval  $J_{n-1} $ in  $ \JJ_{n-1}$  and satisfy \eqref{cond_jn}.
We define $\JJ_n$ by induction on $n$.

For $n=0$, we trivially have that~\eqref{cond_jn} is satisfied
for any interval $J_0 \subset I$.  The point is that
$H(\Delta)\ge 1$ and so there are no intervals $\Delta(L)$  satisfying the height condition $H(\Delta)< 1$. So take $\JJ_0:=\{J_0\}$. For the same
reason~\eqref{cond_jn} with $n=1$ is trivially satisfied for any
interval $J_1$ obtained by subdividing $J_0$ into $R$ closed
intervals of equal length $c_1R^{-1}$. Denote by $\JJ_1$ the
resulting collection of intervals $J_1$.

In general, given $\JJ_n$ satisfying  \eqref{cond_jn} we wish to
construct a nested collection $\JJ_{n+1}$ of intervals $J_{n+1}$ for
which (\ref{cond_jn}) is satisfied with $n$ replaced by $n+1$. By
definition, any interval $J_n$ in $\JJ_n$ avoids intervals
$\Delta(L)$ arising from lines $L$ with height $ H(\Delta) $ bounded above by
$R^{n-1}$.  Since any `new' interval $J_{n+1}$ is to be nested in
some $J_n$,  it is enough to show that $J_{n+1}$ avoids intervals
$\Delta(L)$ arising from lines $L$ with height $ H(\Delta) $  satisfying \eqref{classmore}; that is
$$
R^{n-1}\le H(\Delta)<R^n  \   .
$$
The collection of intervals $\Delta(L)  \in  \RRR$  satisfying this height condition  is precisely  the class $ C(n)\cup C^*(n)$  introduced at the beginning of \S\ref{sec_split}. In other words, it the  precisely the collection $ C(n)\cup C^*(n)$ of  intervals  that come into play  when attempting to construct   $\JJ_{n+1}$ from $\JJ_{n}$. We now proceed
with the construction.

Assume that $n\ge 1$. We subdivide each $J_n$ in $\JJ_n$  into $R$ closed intervals $I_{n+1}$ of equal length $ c_1
R^{-(n+1)}$ and denote by $\II_{n+1}$ the collection of
such intervals. Thus,
$$|I_{n+1}|=c_1 R^{-(n+1)}  \qquad {\rm and }  \qquad \# \II_{n+1}  \, =  \,  R  \, \times \, \# \JJ_{n}   \, .   $$
It is obvious that the construction of $\II_{n+1}$ corresponds
to the splitting procedure associated with the  construction of a  $({\rm I},\vR,\vr)$ Cantor set.

In view of the nested requirement, the collection $\JJ_{n+1}$ which
we are attempting to construct will be a  sub-collection  of
$\II_{n+1}$. In other words, the intervals $I_{n+1}$ represent
possible candidates for $J_{n+1}$. The goal now is simple --- it is
to remove those `bad' intervals $I_{n+1}$ from $\II_{n+1}$  for
which
\begin{equation}  
I_{n+1} \,  \cap \, \Delta(L) \, \neq  \, \emptyset   \ \ \mbox{
for some \ } \Delta(L)   \in C(n)\cup C^*(n)  \ .
\end{equation}
The sought after collection $ \JJ_{n+1}$ consists precisely  of those intervals that survive.   Formally, for $n \ge 1 $  we let
$$
\JJ_{n+1}:=\{I_{n+1}\in\II_{n+1}\;:\;I_{n+1} \,  \cap \,
\Delta(L)=\emptyset\ \mbox{ for any }\Delta(L)\in C(n)\cup C^*(n)\}.
$$



\noindent We claim that these collections of surviving intervals  satisfy the following key statement. It  implies that the  act of removing `bad' intervals from $\II_{n+1}$  is exactly in keeping with the
removal procedure associated with the construction of a  $(J_0,\vR,\vr)$ Cantor set with $\vR $  and $ \vr$  as described in \S\ref{subsec3_2}.

\begin{proposition}\label{prop1}
Let $\epsilon := (ij)^2/ 8 $ and with reference to \S\ref{sec_split}  let
$$
C(n,l):=\bigcup_{k=0}^{\log_2 \! R}C(n,k,l) \ , \hspace*{8ex}
C_1(n):=\bigcup_{k=0}^{\log_2 \! R}C_1(n,k) \ ,
$$

$$
C_2(n):=\bigcup_{k=0}^{\log_2 \! R}C_2(n,k)\quad\mbox{and}\quad
\widetilde{C}_3(n,u):=\bigcup_{k=0}^{\log_2 \! R} \ \bigcup_{v=0}^{\lambda \log_2 \!
R}C_3(n,k,u,v).
$$
Then,  for $R>R_0(\epsilon)$ large enough the following four statements are valid.

\begin{enumerate}
\item For any fixed interval $J_{n-l}\in \JJ_{n-l}$, the intervals  from class $C(n,l) $ with $ n/\lambda\ge l\ge 1$ intersect no more than
$R^{1-\epsilon}$ intervals  $I_{n+1}\in \II_{n+1}  $ with $ I_{n+1}\subset
J_{n-l}$.
\item For any $n\ge 3n_0$ where $n_o$ is defined by \eqref{def_n0} and any fixed  interval $J_{n-n_0}\in
\JJ_{n-n_0}$, the intervals from class $C^*(n)$ intersect  no more
than $R^{1-\epsilon}$ intervals  $I_{n+1}\in \II_{n+1} $  with $
I_{n+1}\subset J_{n-n_0}$.
\item For any fixed interval  $J_n\in \JJ_n$, the  intervals from class  $C_1(n)$ or $C_2(n)$ intersect  no more than
$R^{1-\epsilon}$ intervals  $I_{n+1}\in \II_{n+1}$ with $ I_{n+1}\subset
J_n$.
\item For any fixed interval  $J_{n-u}\in \JJ_{n-u}$, the intervals from class $\widetilde{C}_3(n,u)$  intersect  no more than
$R^{1-\epsilon}$ intervals  $I_{n+1}\in \II_{n+1} $ with $I_{n+1}\subset
J_{n-u}$.
\end{enumerate}
\end{proposition}

%

\medskip

\noindent{\em Remark 1. \,} Note that in Part 1  we have that
$l<n/\lambda$ and in Part 2 we have that  $u$ is bounded above by
\eqref{class3_lbound}. So in either part  we have  that $l,u\le n$ for all
positive values $n$. Therefore the collections  $\JJ_{n-l}$ and
$\JJ_{n-u}$ are well defined.

\medskip

\noindent{\em Remark 2. \,} By definition,  a planar curve $\CCC :=
\CCC_f$ is $C^{(2)}$ non-degenerate if $f \in C^{(2)}(I) $ and there
exits at least one point $ x \in I $ such that
$ f''(x) \neq 0 $.  It will be apparent during the course of establishing Proposition  \ref{prop1}
that the condition on the curvature  is only required when considering Part 2.  For the other parts only the
two times continuously differentiable condition is required.  Thus, Parts 1, 3 and 4  of the proposition remain valid even when the curve is a line. The upshot is that Proposition  \ref{prop1} remains valid for any $C^{(2)} $ curve for which $V_L$ is not too small and for such curves we are  able to establish  the analogue of Theorem \ref{thmfinite}.  We will use this observation when proving  Theorem \ref{thmnvlines}.

\subsection{Dealing with  $\bad(0,1)\cap\CCC$ \label{pain}}

The  construction of the collections $\JJ_{n}$ satisfying Proposition~\ref{prop1} requires that $i>0$.
However, by  making  use of  the fact that $\bad(0,1)\cap \CCC
= (\RR\times \bad) \cap \CCC$,    the case $(i,j)=(0,1)$ can be easily dealt with.
\medskip

\noindent Let $R \ge 2$ be a large integer, and let
\begin{equation}\label{defc10}
c_1 :=\frac{2cR^2}{c_0}  \quad \mbox{ \ where \  } \quad 0 < c < \frac{1}{2R^2}   \; .
\end{equation}

\noindent For a given rational number $p/q$ $(q \ge 1)$,   let $\Delta_\CCC (p/q )$  be the
``interval'' on $\CCC$ defined by
$$
\Delta_\CCC (p/q ):=  \left[f^{-1} \left(\frac{p}{q}\pm
\frac{c}{H(p/q)}\right)
\right]\qquad\mbox{where}\qquad
    H(p/q)  := q^2   \, .
$$
In view of~\eqref{bound_curvat} the inverse function $f^{-1}$ is well
defined. Next observe that the orthogonal  projection of $\Delta_\CCC(p/q)$ onto the $x$-axis is contained in the
interval $\Delta (p/q)$ centered at the point  $f^{-1}(p/q)$ with length
$$
\left|\Delta (p/q)\right| := \frac{2c}{c_0H(p/q)}   \, .
$$

By analogy with \S\ref{fat} the set $\bad^f_c(0,1)$ can
be described as the set of $x \in I$ such that $ x \notin \Delta(p/q)$ for all rationals $ p/q$.   For the sake of consistency
with the $i > 0 $ situation, for $n\ge 0$ let

$$
\cC(n):=\left\{\Delta(p/q)\;:\;p/q\in\QQ  \ \ {\rm and} \  \  R^{n-1}\le
H(p/q)< R^n\right\} .
$$

\noindent Since $\cC(n)=\emptyset $ for $n=0$, the
following analogue of Proposition~\ref{prop1} allows us to deal with
the $i=0$ case. \emph{For $R\ge
4$ and any interval $J_n\in \JJ_n$,  we have that}
\begin{equation}\label{defc10db}
\#\{I_{n+1}\in \II_{n+1}  \, :  I_{n+1}\subset J_n \ {\rm and } \   \Delta(p/q)
\cap I_{n+1}\neq\emptyset  {\rm \  \ for  \ some \   }   \
\Delta(p/q) \in \cC(n) \} \;\le\;
 3    \ .
\end{equation}
In short, it allows us to construct a $(J_0,\vR,\vr)$ Cantor subset   of $\bad^f_c(0,1)$ with
$$\vR:=(R_n)=(R,R,R,\ldots)$$
and
$$\vr:=(r_{m,n})=\left\{\begin{array}{l}3 \quad\mbox{if }m=n;\\[2ex]
0 \quad\mbox{if }m<n.
\end{array}\right.
$$
To establish \eqref{defc10db}  we proceed as follows.  First note that
in view of \eqref{defc10}, we have that
$$
\frac{|\Delta(p/q)|}{|I_{n+1}|}\le 1  \ .
$$
Thus, any single interval  $\Delta(p/q)$ removes at most three
intervals $I_{n+1}$ from $\II_{n+1}$. Next, for any two rationals
$p_1/q_1,p_2/q_2\in \cC(n)$ we have that
$$
\left|f^{-1}\left(\frac{p_1}{q_1}\right)-f^{-1}\left(\frac{p_2}{q_2}\right)\right|
\, \ge \, \frac{1}{|f'(\xi)|q_1q_2}  \, \ge  \,   \frac{1}{c_0} R^{-n}    \,
> \, c_1R^{-n} \, := |J_n|
$$
where $\xi$ is some number between $p_1/q_1$ and $p_2/q_2$. Thus,
there is at most one interval $\Delta(p/q)$ that can possibly
intersect any given interval $J_n$ from $\JJ_n$. This together with
the previous fact establishes \eqref{defc10db}.

\section{Forcing lines to intersect at one point \label{rty} }

From this point onwards, all our effort is geared towards
establishing Proposition~\ref{prop1}. Fix  a generic  interval
$J\subset I$ of length $c_1' R^{-n}$.   Note that the position of
$J$ is not specified and sometimes it may be more illuminating to
picture   $J$ as an interval on $\CCC$.  Consider all intervals
$\Delta(L)$  from the same class (either $C(n,k,l,m)$, $C^*(n,k)$,
$C_1(n,k)\cap C(n,k,0,m)$, $C_2(n,k)$ or $C_3(n,k,u,v)$) with
$\Delta(L)\cap J\neq \emptyset$. The overall aim of this section is
to determine conditions on the size of $c'_1$ so that the associated
lines $L$ necessarily intersect at single  point.

\subsection{Preliminaries: estimates for $F_L$ and $F_L'$}

Let
\begin{equation}\label{lbound_c1}
c_1'\ge 2Kc^{1/2}\cdot 2^{-k}R.
\end{equation}
This condition guarantees that any interval $\Delta(L)\in C(n,k,l)$
(or $\Delta(L)\in C^*(n,k)$) has length smaller than $|J|$. Indeed,
$$
|\Delta(L)| =   2Kc^{1/2}\cdot (H(\Delta))^{-1}
\stackrel{\eqref{class_prop1}}\le 4Kc^{1/2}R\cdot 2^{-k}R^{-n} \le
|J| \, .
$$

\noindent In this section we obtain various  estimates for $|F_L(x)|$ and
$|F'_L(x)|$ that are valid for any  $x\in J$.   Recall,  $x_0$ is as usual the point at which $|F_L'(x)|$ attains its minimum with $x\in \Delta(L)$.

\begin{lemma}\label{lem_vp} Let $ 0   \leq m \leq \lambda\log_2 R$, $l \ge 0$ and
 $c_1'$ be a positive parameter such that
\begin{equation}\label{cond_c1d}
8\vC c_1'R^{-n}\le 2^{-m} R^{-\lambda l}.
\end{equation}
Let $J\subset I$ be an interval of length $c_1'R^{-n}$. Let
$\Delta(L)$ be any interval from class $C(n,k,l,m)$ such that
$\Delta(L)\cap J\neq \emptyset$. Then for any $x\in J$ we have
$|F'_L(x)|\asymp V_L$ and
\begin{equation}\label{bound_f}
|F_L(x)|\le 5|J|V_L.
\end{equation}

\end{lemma}
\proof A consequence of Taylor's formula is that
\begin{eqnarray}\label{ineq_lem2}
|F'_L(x)-V_L| &=& |A-Bf'(x)-V_L|=|x-x_0|\cdot
|-Bf''(\tilde{x})|\nonumber\\[2ex]
&\le&
(c_1'+2Kc^{1/2}R)R^{-n}\cdot \vC\max\{|A|,|B|\}\nonumber\\
&\stackrel{\eqref{lbound_c1}}\le&
2c_1'R^{-n}\cdot \vC\max\{|A|,|B|\}
\end{eqnarray}
where $\tilde{x}$ is some point between $x$ and $x_0$. Then by~\eqref{lbound_c1} and~\eqref{cond_c1d} together with the fact that $\Delta(L)\in C(n,k,l,m)$ we get that
$$
|F'_L(x)-V_L|\le \frac12\cdot 2^{-m-1}R^{-\lambda l}\max\{|A|,|B|\}\stackrel{\eqref{class0m_prop}}{\le} \frac12 V_L  \, .
$$
In other words, $|F'_L(x)|\asymp V_L$. Then
$$
|F_L(x)|\le |F_L(x_1)|+ |x-x_1|\cdot|F_L'(\tilde{x})|\le
\frac{c}{\maxab} + 4|J|V_L
$$
where $x_1$ is the center of $\Delta(L)$ and $\tilde{x}$ is some
point between $x$ and $x_1$. However
$$  c  \, ({\maxab})^{-1} =    c^{1/2} V_L (H(\Delta))^{-1} \stackrel{\eqref{class_prop1}}\le c^{1/2}R\cdot
R^{-n} V_L \le |J|V_L
$$
and as a consequence,~\eqref{bound_f} follows.
\endproof

\begin{lemma}\label{lem_vp2}
Assume $c_1'$ does not satisfy \eqref{cond_c1d}.
Let $J\subset I$ be an interval of length $c_1'R^{-n}$. Let
$\Delta(L(A,B,C))\in C(n,k,l,m)$ such that
$\Delta(L)\cap J\neq \emptyset$. Then for any $x\in J$ we have
\begin{equation}\label{bound_fp}
|F_L(x)|\le 30 \vC |J|^2 \max\{|A|,|B|\}
\end{equation}
and
\begin{equation}\label{bound_dfp}
|F'_L(x)|\le 10\vC |J| \max\{|A|,|B|\}.
\end{equation}
\end{lemma}

\proof In view of \eqref{ineq_lem2} it follows that
$$
|F'_L(x)|\le 2c_1'R^{-n}\vC \maxhab + V_L  \, .
$$
By \eqref{class0m_prop} we have that
$$
V_L\le 2^{-m}R^{-\lambda l}\maxhab \le 8\vC |J|\maxhab.
$$
Combining these estimates gives~\eqref{bound_dfp}.

To establish  inequality~\eqref{bound_fp} we use Taylor's formula.  The latter implies the existence of  some point  $\tilde{x}$ between $x$ and $x_1$ such that
\begin{eqnarray*}
|F_L(x)|  & \le  &  |F_L(x_1)|+|x-x_1||F'_L(x_1)|+\frac12|x-x_1|^2
|-Bf''(\tilde{x})|  \\
& \stackrel{\eqref{bound_dfp}}\le & \frac{c}{\maxab}+20\vC |J|^2
\max\{|A|,|B|\} +2 \vC |J|^2 \max\{|A|,|B|\}.
\end{eqnarray*}
This together with the fact that the  first of the three terms on the r.h.s. is bounded above by $c^{1/2} V_L
(H(\Delta))^{-1} \le 8\vC |J|^2 \maxhab$   yields~\eqref{bound_fp}.
\endproof

\medskip

The next lemma provides an estimate for $F_L(x)$ and $F'_L(x)$ in
case $\Delta(L)$ is of Type 2.

\begin{lemma}\label{lem_tp2vp}
Let $c_1'$ be a positive parameter such that
\begin{equation}\label{cond_tpc1d}
1\le \vC c_1'\quad\mbox{and}\quad R^2c\le \vC c_1^{'2}.
\end{equation}
Let $J\subset I$ be an interval of length $c_1'R^{-n}$. Let
$\Delta(L)$ be any interval from class $C^*(n,k)$ such that
$\Delta(L)\cap J\neq \emptyset$. Then for any $x\in J$ we have
\begin{equation}\label{bound_tp2fp}
|F_L(x)|\le 9 \vC |J|^2 \max\{|A|,|B|\}
\end{equation}
and
\begin{equation}\label{bound_tp2dfp}
|F'_L(x)|\le 3\vC |J| \max\{|A|,|B|\}.
\end{equation}
\end{lemma}

\proof As in the previous two lemmas a simple consequence of
Taylor's formula is that there exists $\tilde{x}$
between $x$ and $x_0$ such that:
$$
|F_L'(x)|\le V_L + |x-x_0|\cdot |-Bf''(\tilde{x})|
\stackrel{\eqref{type2_eq}}\le R^{-n}\maxhab + 2\vC|J| \maxhab
$$
which by~\eqref{cond_tpc1d} leads to~\eqref{bound_tp2dfp}. For the
first inequality, by Taylor's formula  we have that
\begin{equation} \label{ju}
|F_L(x)| \le \frac{c}{\maxab} + 8\vC|J|^2 \maxhab
\end{equation}
On the other hand by~\eqref{class_prop1} we have that
$$
H(\Delta) = (\maxab |B|)^{1/2}\ge R^{n-1}
$$ and so
$$
\frac{c}{\maxab} \ \le \  \frac{R^2c |B|}{R^{2n}}
 \ \stackrel{\eqref{cond_tpc1d}}\le \  \vC |J|^2 \maxhab  \, .
$$

 \noindent This together with \eqref{ju} yields~\eqref{bound_tp2fp} . \endproof

\subsection{Avoiding Parallel lines}\label{subseq1}

Consider all lines $L_1,L_2,\cdots$ such that the corresponding
intervals $\Delta(L_1),\Delta(L_2),\cdots$ belong to the same class
and intersect  $J$.  Recall,  $|J|:=  c_1' R^{-n}$.   In this section,  we determine conditions on
$c_1'$ which ensure that none of the lines $L_i$  are parallel to one another.

\medskip

\noindent{\em Remark 1. \,}  For the sake of clarity and  to minimize  notation,  throughout the rest of the paper we will often write $V_1$, $V_2,
\cdots$ instead of $V_{L_1}$, $V_{L_2},\cdots$ when there is no risk of ambiguity.

\medskip

\begin{lemma} \label{5}
Assume that there are at least two parallel lines $L_1(A_1,B_1,C_1),
L_2(A_2,B_2,C_2)$ such that $\Delta(L_1)\cap J\neq\emptyset$ and
$\Delta(L_2)\cap J\neq\emptyset$. If $\Delta(L_1),\Delta(L_2)\in C(n,k,l,m)$
and~\eqref{cond_c1d} is satisfied then
\begin{equation}\label{paral_cond}
c'_1V_1\min\{|A_1|,|B_1|\}\gg R^n.
\end{equation}
If $\Delta(L_1),\Delta(L_2)\in C(n,k,l,m)$ and~\eqref{cond_c1d} is false
or $\Delta(L_1),\Delta(L_2)\in C^*(n,k)$ and~\eqref{cond_tpc1d} is true then
\begin{equation}\label{paral_cond2}
c'_1\sqrt{|A_1||B_1|}\gg R^n.
\end{equation}

\end{lemma}

\proof Assuming  that $L(A_1,B_1,C_1),
L(A_2,B_2,C_2)$ are parallel  implies that  $A_2=t A_1, B_2=t B_1,
t\in\QQ$. Without loss of generality,  assume that $|t|\le 1$. This implies that $|A_1|\ge |A_2|$ and $|B_1|\ge |B_2|$.
Then for an arbitrary point $x\in J$, we have
\begin{equation} \label{yud}
|tC_1-C_2| = |tF_{L_1}(x) - F_{L_2}(x)|  \, .
\end{equation}
The denominator of $t$ divides both  $A_1$ and $B_1$ so $t$ is at
most $\min(|A_1|,|B_1|)$. Therefore the l.h.s. of \eqref{yud} is at
least $(\min\{|A_1|,|B_1|\})^{-1}$.

If $c_1'$ satisfies~\eqref{cond_c1d} then the conditions of
Lemma~\ref{lem_vp} are true. Therefore $V_1\asymp V_2$ and r.h.s. of \eqref{yud} is
at most $5|J|(V_1+V_2)\ll c_1'V_1 R^{-n}$. This together with the
previous estimate  for the l.h.s. of \eqref{yud} gives~\eqref{paral_cond}.  To  establish the remaining part of the lemma, we exploit either Lemma~\ref{lem_vp2} or Lemma~\ref{lem_tp2vp}
to show that
$$
{\rm r.h.s. \ of \ \eqref{yud}} \ll |J|^2\max\{|A_1|,|B_1|\} =
(c_1'R^{-n})^2\max\{|A_1|,|B_1|\}.
$$
This together with the
previous estimate for the l.h.s. of \eqref{yud}  gives~\eqref{paral_cond2}.
\endproof

\vspace*{2ex}

The upshot of Lemma \ref{5}  is that there are no  parallel lines in the same class passing through a generic   $J$ of length $c_1' R^{-n}$ if  $c_1'$ is  chosen to be sufficiently small so that \eqref{paral_cond} and \eqref{paral_cond2} are violated.

\subsection{Ensuring lines  intersect at one point}\label{sec_inter}

Recall, our aim is to determine conditions on
$c_1'$ which ensure that all lines $L$ associated with intervals $\Delta(L)$  from the same class with   $\Delta(L)\cap J\neq \emptyset$ intersect at one point. We will use the following well-known fact. For $i=1,2,3$, let  $L_i(A_i,B_i,C_i)$ be a line given  by the  equation $A_ix-B_iy+C_i=0$.
The lines do not intersect at a single  point if and only if
$$
\det\left(\begin{array}{ccc} A_1&B_1&C_1\\
A_2&B_2&C_2\\
A_3&B_3&C_3
\end{array}\right)\neq 0.
$$

\noindent Suppose there are at  least three intervals
$\Delta(L_1),\Delta(L_2),\Delta(L_3)$ from the same class (either
$C(n,k,l)$, $C_1(n,k), C_2(n,k), C_3(n,k,u,v)$ or $C^*(n,k)$) that intersect $J$ but
 the corresponding lines $L_1,L_2$ and $L_3$ do not intersect at a single
point. Then
$$
\left|\det\left(\begin{array}{ccc} A_1&B_1&C_1\\
A_2&B_2&C_2\\
A_3&B_3&C_3
\end{array}\right)\right|\ge 1.
$$

\noindent Choose an arbitrary point $x\in J$. Firstly assume that $J$
satisfies~\eqref{cond_c1d} and that the intervals
$\Delta(L_1),\Delta(L_2),\Delta(L_3)$ are of Type 1. Then
Lemma~\ref{lem_vp} implies that
$$
F_{L_1}(x)\ll |J|V_1.
$$
The same inequalities are true for  $F_{L_2}(x)$ and
$F_{L_3}(x)$. So we obtain that
$$
\left(\begin{array}{ccc} A_1&B_1&C_1\\
A_2&B_2&C_2\\
A_3&B_3&C_3
\end{array}\right)\cdot
\left(\begin{array}{c}x\\f(x)\\1\end{array}\right)\ll
\left(\begin{array}{c}|J|V_1\\
|J|V_2\\
|J|V_3\end{array}\right).
$$
By applying Cramer's rule  to the third row, we find that
$$
|J| \; (|V_1(A_2B_3-A_3B_2)|+|V_2(A_1B_3-A_3B_1)|+|V_3(A_1B_2-A_2B_1)|)\gg
1.
$$
Without loss of generality assume that the first term on the l.h.s. of this
inequality is the largest of the three terms. Then
\begin{equation}\label{triang}
c_1'|V_1(A_2B_3-A_3B_2)|\gg R^n.
\end{equation}
In other words, if the lines $L_1,L_2$ and $L_3$ do not intersect at
one point and~\eqref{cond_c1d} is true for a given  $c_1'$
then~\eqref{triang} must also hold.

\medskip

If~\eqref{cond_c1d} is not true or the  intervals $\Delta(L_1),
\Delta(L_2),\Delta(L_3)$ are of Type 2 then we apply either
Lemma~\ref{lem_vp2} or Lemma~\ref{lem_tp2vp}. Together with Cramer's
rule,  we obtain that
$$
|J|^2(\max\{|A_1|,|B_1|\}|A_2B_3-A_3B_2|+\max\{|A_2|,|B_2|\}|A_1B_3-A_3B_1|
$$$$
+\max\{|A_3|,|B_3|\}|A_1B_2-A_2B_1|)\gg 1.
$$
Without loss of generality assume that the first of the three terms  on the l.h.s. of this inequality is the largest. Then,  we obtain that
\begin{equation}\label{triang2}
c_1'\sqrt{|\max\{|A_1|,|B_1|\}(A_2B_3-A_3B_2)|}\gg R^n.
\end{equation}

We now investigate the ramifications of the conditions~\eqref{triang} and~\eqref{triang2} on specific classes of intervals.

\subsubsection{Case $\Delta(L_1),\Delta(L_2),\Delta(L_3)\in C(n,k,l,m), l\ge 1$}

We start by estimating the difference between $\frac{A_1}{B_1}$ and
$\frac{A_2}{B_2}$. By \eqref{class_prop2} we have  that
\begin{equation}\label{class_eq4}
\left|\frac{A_1}{B_1} - \frac{A_2}{B_2}\right|\le
\left|\frac{A_1}{B_1} -
f'(x_{01})\right|+|f'(x_{01})-f'(x_{02})|+\left|f'(x_{02}) -
\frac{A_2}{B_2}\right|\ll R^{-\lambda l} + |J|
\end{equation}
where $x_{01}$ and $x_{02}$ are given by $V_1 := |A_1-B_1f'(x_{01})|$ and $V_2 := |A_2-B_2f'(x_{02})|$ respectively.


\vspace*{2ex}

$\bullet$ {\it Assume that~\eqref{cond_c1d} is satisfied}. This means
that $|J|\ll R^{-\lambda l}$. We rewrite~\eqref{triang} as
$$
c_1'|V_1B_2B_3|\left|\frac{A_2}{B_2}-\frac{A_3}{B_3}\right|\gg R^n.
$$
Then in view of \eqref{class_eq2}, \eqref{class_prop2} and
\eqref{class_eq4} it follows that

\begin{eqnarray*}
R^{n}  & \ll & c'_1R^{-\lambda
l}\left(\frac{2^kc^{\frac12}}{R}R^{n+\lambda
(l+1)}\right)^{\frac{3i}{i+1}}\cdot R^{-\lambda l}  \\[2ex] & \stackrel{\eqref{cond_c1d}}\ll &  c'_1R^{n-\frac{j-i}{i+1}n}\cdot
R^{-\frac{2-i}{i+1}\lambda l}\cdot
\left(\frac{2^kc^{\frac12}}{R}R^\lambda\right)^{\frac{3i}{i+1}} \, .
\end{eqnarray*}
Since by assumption $i\le j$, the last inequality implies that if~\eqref{triang} holds then
$$
c'_1\gg
R^{l\lambda\frac{2-i}{i+1}}\cdot\left(\frac{2^kc^{\frac12}}{R}R^\lambda\right)^{-\frac{3i}{i+1}}.
$$
Hence,  the condition
$$
c'_1\ll
R^{l\lambda\frac{2-i}{i+1}}\cdot\left(\frac{2^kc^{\frac12}}{R}R^\lambda\right)^{-\frac{3i}{i+1}}
$$
will contradict the previous inequality and imply that~\eqref{triang} is not satisfied.
Note that similar  arguments imply  that if~\eqref{paral_cond} holds then
$$
R^n\ll c_1'R^{-\lambda l}\left(\frac{2^kc^{\frac12}}{R}R^{n+\lambda
(l+1)}\right)^{\frac{2i}{i+1}}=c_1'R^{n-\frac{j}{i+1}n}R^{-\frac{j}{i+1}\lambda
l}\cdot
\left(\frac{2^kc^{\frac12}}{R}R^\lambda\right)^{\frac{2i}{i+1}}   \, .
$$
It follows that the condition
$$
c'_1\ll
R^{l\lambda\frac{j}{i+1}}\cdot\left(\frac{2^kc^{\frac12}}{R}R^\lambda\right)^{-\frac{3i}{i+1}}.
$$
will contradict the previous inequality and imply that~\eqref{paral_cond} is not satisfied.

The upshot is that for $\lambda$ satisfying~\eqref{cond_lambda} the following
condition on $c_1'$
\begin{equation}\label{cond_c1_2}
c_1'\le \delta\cdot R^l\cdot
\left(\frac{2^kc^{\frac12}}{R}R^\lambda\right)^{-\frac{3i}{i+1}}
\end{equation}
will contradict both~\eqref{triang} and~\eqref{paral_cond}.  Here $\delta=\delta(i,j,c_0,C_0)>0$ is the  absolute unspecified constant within  the previous inequalities involving the  Vinogradov symbols.  In other words,
if $c_1'$ satisfies \eqref{cond_c1_2},  then  the lines $L_i$ associated with the  intervals
$\Delta(L_i)\in C(n,k,l,m)$ with $ l\ge 1$ such that   $\Delta(L_i)  \cap J \neq \emptyset$     intersect at a single  point.

\vspace*{2ex}

$\bullet$ {\it Assume that~\eqref{cond_c1d} is false}. In this case
$R^{-\lambda l}\ll R^\lambda |J|$. In view of~\eqref{class_eq1} we
have that $|A_1|\asymp|B_1|$ and inequality~\eqref{triang2} implies  that
$$
c_1'\sqrt{|B_1B_2B_3|\left|\frac{A_2}{B_2}-\frac{A_3}{B_3}\right|}\gg
R^n.
$$
In view of \eqref{class_eq2} and \eqref{class_eq4},  it follows that
to
$$
R^{n}\ll (c'_1)^{\frac32}R^{\lambda/2 - n/2}
\left(\frac{2^kc^{\frac12}}{R}R^{n+\lambda
(l+1)}\right)^{\frac{3i}{2(i+1)}}
$$
which is equivalent to
$$
R^n \ll
c'_1R^{\lambda/3}\left(\frac{2^kc^{\frac12}}{R}R^\lambda\right)^{\frac{i}{(i+1)}}
(R^{n+\lambda l})^{\frac{i}{i+1}}.
$$
This together with that fact that  $i\le 1/2$ and $\lambda l\le n$ implies that
$$
c_1'\gg R^{\frac{j}{i+1}\lambda
l}\cdot\left(\frac{2^kc^{\frac12}}{R}R^\lambda\right)^{-\frac{i}{i+1}}R^{-\frac{\lambda}{3}}.
$$
By similar arguments,  estimate~\eqref{paral_cond2} implies that
$$
c_1' \gg R^{\frac{j}{i+1}\lambda
l}\cdot\left(\frac{2^kc^{\frac12}}{R}R^\lambda\right)^{-\frac{i}{i+1}}.
$$
The upshot is that for $\lambda$ satisfying~\eqref{cond_lambda}, we obtain a contradiction to both these upper bound inequalities  for $c_1'$  and thus to~\eqref{triang2} and~\eqref{paral_cond2},   if
\begin{equation}\label{cond_c1_22}
c_1'\le \delta\cdot R^l\cdot
\left(\frac{2^kc^{\frac12}}{R}R^\lambda\right)^{-\frac{i}{i+1}}
R^{-\frac{\lambda}{3}}.
\end{equation}
In other words,  if $c_1'$ satisfies \eqref{cond_c1_22} but
not~\eqref{cond_c1d},  then the lines $L_i$ associated with the intervals
 $\Delta(L_i)\in
C(n,k,l,m)$ with  $l\ge 1$ such that $\Delta(L_i)  \cap J \neq \emptyset$ intersect at a single point.

\subsubsection{Case $\Delta(L_1),\Delta(L_2),\Delta(L_3)\in C^*(n,k)$}

For this class of intervals we will  eventually make use of  Lemma \ref{lem_tp2vp}.  With this in mind, \emph{we assume that~\eqref{cond_tpc1d} is valid.}  A consequence  of~\eqref{cond_tpc1d} is that  $R^{-n}\ll |J|$. It is readily verified that in the case under
consideration, the analogy to~\eqref{class_eq4} is given by
$$
\left|\frac{A_1}{B_1} - \frac{A_2}{B_2}\right| \, \ll  \,
\frac{V_1}{B_1}+|J|+\frac{V_2}{B_2}\stackrel{\eqref{type2_eq}}\ll  \,
|J|.
$$
Then by using~\eqref{class_eq2_2}, we find that
inequality~\eqref{triang2} implies that
$$
c_1' \gg R^{\frac{j}{1+i}n}\cdot
\left(\frac{2^k}{R}\right)^{-\frac{2i}{i+1}}.
$$
Similarly,  inequality~\eqref{paral_cond2} implies the same upper bound for $c_1'$.
Thus if $c_1'$ satisfies the condition
\begin{equation}\label{cond_c1_tp2}
c_1'\le \delta\cdot R^{\frac{j}{1+i}n}\cdot
\left(\frac{2^k}{R}\right)^{-\frac{2i}{(i+1)}}  \ ,
\end{equation}
 we obtain a contradiction to both~\eqref{triang2}
and~\eqref{paral_cond2}.

\vspace*{3ex}


\subsubsection{Case $\Delta(L_1),\Delta(L_2),\Delta(L_3)\in C_1(n,k)\cap C(n,k,0,m)$ \label{need}}
In view of  \eqref{class0_aeq} and \eqref{class1},   inequality
\eqref{triang} implies  that
$$
R^n\ll
c_1'\left(\frac{2^kc^{\frac12}}{R}R^{n+\lambda}\right)^{\frac{3i}{i+1}}\;\;\stackrel{i\le
1/2}\ll\;\;c_1'\left(\frac{2^kc^{\frac12}}{R}R^{n+\lambda}\right)   \, .
$$
Hence, if $c_1'$ satisfies the condition
\begin{equation}\label{cond_c1_3}
c'_1\le \delta\cdot
\left(\frac{2^kc^{\frac12}}{R}R^\lambda\right)^{-1}
\end{equation}
we obtain a contradiction to~\eqref{triang}. Note that the same
upper bound inequality  for $c_1'$  will also contradict
\eqref{paral_cond}.

For the  class $C_1(n,k)$ as well as all other subclasses of
$C(n,k,0)$,   when consider the intersection with a generic
interval $J$  of length $c_1' R^{-n}$  the constant  $c_1'$  will
always satisfy~\eqref{cond_c1d}.   \emph{Therefore, without loss of
generality we assume that  $c_1'$ satisfies~\eqref{cond_c1d}}.

\subsubsection{Case $\Delta(L_1),\Delta(L_2),\Delta(L_3)\in C_2(n,k)$}
By \eqref{class0_beq}, \eqref{class2_aeq} and
\eqref{class2_veq},   inequality \eqref{triang} implies that
$$
R^{n}\ll
c_1'\left(\frac{2^kc^{\frac12}}{R}R^{n+\lambda}\right)^{\frac{i}{j+1}+\frac{2j}{j+1}}=c_1'\left(\frac{2^kc^{\frac12}}{R}R^{n+\lambda}\right).
$$
It is now easily verified that if  $c'_1$  satisfies inequality  \eqref{cond_c1_3} then we obtain a contradiction to~\eqref{triang}.

\subsubsection{Case $\Delta(L_1),\Delta(L_2),\Delta(L_3)\in C_3(n,k,u,v)$}
By \eqref{class3_beq} and \eqref{class3_aeq},  inequality
\eqref{triang} implies that
$$
R^{n}\ll
c_1'\left(\frac{2^kc^{\frac12}}{R}R^{n}\right)^{\frac{i}{j+1}+\frac{2j}{j+1}}\!\!\!\!\cdot
R^{-\frac{2\lambda uj}{j+1}+\frac{\lambda uij}{j+1}}R^{\lambda i}
\ll c_1'\left(\frac{2^kc^{\frac12}}{R}R^n\right) R^{\lambda
i-\lambda uj}.
$$
Hence, if $c_1'$ satisfies the condition
$$
c'_1\le \delta\cdot
\left(\frac{2^kc^{\frac12}}{R}\right)^{-1}R^{\lambda uj-\lambda i}
\, ,
$$
we obtain a contradiction to~\eqref{triang}.
It is easily verified that if  $c'_1$  satisfies this lower bound inequality, then  we also obtain a contradiction to \eqref{paral_cond} as well.

It follows by~\eqref{cond_lambda} that  $\lambda\ge 1/j$ and
therefore the above  lower bound inequality for $c_1'$ is true if
\begin{equation}\label{cond_c1_4}
c'_1\le \delta\cdot R^u\cdot \frac{R^{1-\lambda i}}{2^kc^{\frac12}}.
\end{equation}

The upshot of this section is as follows. Assume that $\Delta(L_1),
\Delta(L_2),\Delta(L_3)$ all intersect $J$ and belong to the same
class. Then for each class,   specific conditions for
$c_1'$ have been determined that force the corresponding lines $L_1,L_2$ and $L_3$ to
intersect at  a single  point. These conditions are \eqref{cond_c1d},
\eqref{cond_tpc1d}, \eqref{cond_c1_2}, \eqref{cond_c1_22},
\eqref{cond_c1_tp2}, \eqref{cond_c1_3} and \eqref{cond_c1_4}.

\section{Geometrical properties of  pairs $(A,B)$}

Consider two intervals $\Delta(L_1), \Delta(L_2)\in\RRR$ where
the associated lines $L_1(A_1,B_1,C_1)$ and $L_2(A_2,B_2,C_2)$ are not
parallel. Denote by $P$ the point of intersection $L_1\cap L_2$. To begin with we
investigate the  relationship between
$P,\Delta(L_1)$ and $\Delta(L_2)$.

It is easily seen that
$$
P=\left(\frac{p}{q},\frac{r}{q}\right)=\left(\frac{C_2B_1-C_1B_2}{A_1B_2-A_2B_1},\frac{A_1C_2-A_2C_1}{A_1B_2-A_2B_1}\right);
\quad (p,r,q)=1.
$$
Therefore
\begin{equation}\label{q_eq}
A_1B_2-A_2B_1=tq,\quad  C_1B_2-C_2B_1=-tp,  \quad  A_1C_2-A_2C_1=tr
\end{equation}
for some integer $t$.
Let $x_1$ and $x_2$ be two arbitrary points on $\Delta(L_1)$ and
$\Delta(L_2)$. Since $P\in L_1\cap L_2$,  it  follows  that
$$
\begin{array}{l}
A_1(x_1-\frac{p}{q})-B_1(f(x_1)-\frac{r}{q})=
F_{L_1}(x_1), \\[2ex]
A_2(x_2-\frac{p}{q})-B_2(f(x_2)-\frac{r}{q})= F_{L_2}(x_2)  \, .
\end{array}
$$
By Taylor's formula the second equality can be written as
$$
A_2\left(x_1-\frac{p}{q}\right)-B_2\left(f(x_1)-\frac{r}{q}\right)=F_{L_2}(x_2)+(x_1-x_2)F'_{L_2}(\tilde{x}),
$$
where $\tilde{x}$ is some point between $x_1$ and $x_2$.  This together with the first equality gives
$$
\left(\begin{array}{cc}A_1&-B_1\\
A_2&-B_2
\end{array}\right)\cdot \left(\begin{array}{c}x_1-\frac{p}{q}\\
f(x_1)-\frac{r}{q}
\end{array}\right)=\left(\begin{array}{c}F_{L_1}(x_1)\\
F_{L_2}(x_2) + (x_1-x_2)F'_{L_2}(\tilde{x})
\end{array}\right).
$$
which  on applying  Cramer's rule leads to
\begin{equation}\label{cond100}
x_1-\frac{p}{q}=\frac{B_1(F_{L_2}(x_2)+(x_1-x_2)F'_{L_2}(\tilde{x}))-B_2F_{L_1}(x_1)}{\det\vA}
\end{equation}
and
\begin{equation}\label{cond101}
f(x_1)-\frac{r}{q}=\frac{A_1(F_{L_2}(x_2)+(x_1-x_2)F'_{L_2}(\tilde{x}))-A_2F_{L_1}(x_1)}{\det\vA}.
\end{equation}
Here $$\det\vA:= -A_1B_2+ A_2B_1    \stackrel{\eqref{q_eq}}= -tq \, . $$

Now assume that both intervals $\Delta(L_1)$ and $\Delta(L_2)$
belong to the same class and intersect a fixed generic  interval
$J$ of length $c_1'R^{-n}$. Then, we exploit the fact that $x_1,x_2$ can both be
taken in $J$. Firstly consider the case that $J$ satisfies
\eqref{cond_c1d} and $\Delta(L_1), \Delta(L_2)$ are of Type~1. Then
by Lemma~\ref{lem_vp}
$$ \mbox{ $F'_{L_2}(\tilde{x})\asymp V_2$, \
$F_{L_1}(x_1)\ll |J|V_1$, \  $F_{L_2}(x_2)\ll |J|V_2$ \ and \ $|x_1-x_2|\le
|J|=c_1'R^{-n}$.} $$  This together with  \eqref{cond100} and \eqref{cond101} implies that
\begin{equation}\label{cond_dist}
\begin{array}{rcl}
\displaystyle \frac{|B_1|V_2+|B_2|V_1}{R^n}&\gg& \displaystyle
\frac{|qx_1-p|}{c'_1},\\[4ex]
\displaystyle \frac{|A_1|V_2+|A_2|V_1}{R^n}&\gg& \displaystyle
\frac{|qf(x_1)-r|}{c'_1}.
\end{array}
\end{equation}

\noindent If $J$ does not satisfy~\eqref{cond_c1d} and
$\Delta(L_1),\Delta(L_2)$ are of Type 1  we make use of
Lemma~\ref{lem_vp2} to estimate the  size of  $F_{L_2}(x_2)$, $F'_{L_2}(\tilde{x})$ and
$F_{L_1}(x_1)$. This together with \eqref{cond100} and \eqref{cond101} implies that
\begin{equation}\label{cond_dist2}
\begin{array}{rcl}
\displaystyle
\frac{(|B_1|\max\{|A_2|,|B_2|\}+|B_2|\max\{|A_1|,|B_1|\})}{R^{2n}}&\gg&
\displaystyle \frac{|qx_1-p|}{(c'_1)^2},\\[4ex]
\displaystyle
\frac{(|A_1|\max\{|A_2|,|B_2|\}+|A_2|\max\{|A_1|,|B_1|\})}{R^{2n}}&\gg&
\displaystyle \frac{|qf(x_1)-r|}{(c'_1)^2}.
\end{array}
\end{equation}

\noindent On making use of  Lemma~\ref{lem_tp2vp}, it is easily verified  that the same
inequalities  are valid when  $\Delta(L_1)$, $\Delta(L_2)$ are of Type 2
and  $J$ satisfies~\eqref{cond_tpc1d}.

\subsection{The case $P$ is close to $\CCC$}\label{sec_small_dist}

We consider the situation when the point $P=(p/q,r/q)$ is situated close to the curve $\CCC$. More precisely, assume that there exists at
least one point $(x,f(x))\in\CCC$ such that,
$$
\left|x-\frac{p}{q}\right|<\frac{c}{2}\cdot q^{-1-i},\quad
\left|f(x)-\frac{r}{q}\right|<\frac{c}{2}\cdot q^{-1-j}.
$$
We show that every such point $x$ is situated inside
$\Delta(L_0)$ for some line $L_0$  passing through $P$. Indeed, each
line $L(A,B,C)$ which passes through $P$ will satisfy the equation
$Ap - Br + Cq= 0$. By Minkowski's Theorem there exists an integer
non-zero solution $A_0,B_0,C_0$ to this equation such that
$$
|A_0|<q^i;\quad |B_0|<q^j.
$$
Then
$$
|F_{L_0}(x)| =
|A_0x-B_0f(x)+C_0|=\left|A_0\left(x-\frac{p}{q}\right)-B_0\left(f(x)-\frac{r}{q}\right)\right|\le
cq^{-1}
$$
since $|A_0\cdot\frac{p}{q}-B_0\cdot\frac{r}{q}+C_0|=0$. In other
words, the point $x \in \Delta(L_0)$.

\subsection{The figure $F$ \label{nicefig} }

Consider all intervals $\Delta(L_t(A_t,B_t,C_t))$ from the same
class (either $C(n,k,l,m)$ with $l\ge 1$, $C^*(n,k)$, $C_1(n,k)\cap
C(n,k,0,m),C_2(n,k)$ or $C_3(n,k,u,v)$) which intersect a generic  interval
$J$ of length $c_1'R^{-n}$. In this section we
investigate the implication of this on the coefficients of the corresponding lines  $L_t$.

In \S\ref{rty} we have shown that under  certain conditions on
$c'_1$ all the corresponding lines $L_t$ intersect at one point.
Assume now that the  appropriate conditions are satisfied  -- this
depends of course on the class of intervals under consideration. Let
$P=(p/q,r/q)$ denote the point of intersection of
the lines $L_t$.  Then the triple $(A_t,B_t,C_t)$ will satisfy
the equation
$$
A_tp-B_tr+C_tq=0  \,  \qquad A_t,B_t,C_t\in \ZZ.
$$
Hence the  points  $(A_t,B_t) \in \ZZ^2$ form a lattice $\LL$  with fundamental domain of area equal to~$q$.

\medskip

Let $x_t$ be the point of minimum of $|F'_{L_t}(x)|$ on
$\Delta(L_t)$. Define
$$
\omega_x(P,J):=\max_{t}\{|qx_t-p|\} \quad {\rm and \ } \quad
\omega_y(P,J):=\max_{t}\{|qf(x_t)-r|\}   \ .
$$
Furthermore, let $t_1$ (resp. $t_2$) be    the integer at which the maximum  associated with $ \omega_x $ (resp. $\omega_x$) is attained; i.e.
$$ |qx_{t_1}-p| = \omega_x(P,J)  \quad {\rm and \ } \quad   |qf(x_{t_2})-r| = \omega_y(P,J)  \, .  $$
We  now consider several cases.

\subsubsection{Interval $J$ satisfies~\eqref{cond_c1d} and intervals $\Delta(L_1)$ are of Type 1}

Assume that the  interval $J$ satisfies~\eqref{cond_c1d}. Then on
applying~\eqref{cond_dist} with respect to the pair of intervals
($\Delta(L_t)$, $\Delta(L_{t_1})$) and ($\Delta(L_t)$,
$\Delta(L_{t_2})$), we find that the following   two conditions are satisfied:
\begin{equation}\label{dist_x}
\frac{|B_{t_1}V_t|+|B_tV_{t_1}|}{R^n}\ge
v_x:=\frac{\omega_x(P,J)}{c'_1c_x(\vC,c_0,i,j)} \ \qquad t\neq t_1
\end{equation}
\begin{equation}\label{dist_y}
\frac{|A_{t_2}V_t|+|A_tV_{t_2}|}{R^n}\ge
v_y:=\frac{\omega_y(P,J)}{c'_1c_y(\vC,c_0,i,j)} \ \qquad t\neq t_2 \, ,
\end{equation}
where $c_x(\vC,c_0,i,j)$ and $c_y(\vC,c_0,i,j)$ are
constants dependent only on $\vC,c_0, i$ and $j$.


Firstly consider  inequality~\eqref{dist_x}. Since all intervals
$\Delta(L_t)$ lie in the same class ($C(n,k,l,m)$ with $l\ge 1$,
$C_1(n,k)\cap C_(n,k,0,m), C_2(n,k)$ or $C_3(n,k,u,v)$),  then by
either~\eqref{class2_b12} or~\eqref{asymp_prop} we have
$V_{t_1}\asymp V_t$. This together with~\eqref{class_prop1}
substituted into~\eqref{dist_x} gives
$$
v_x\le
\frac{|B_{t_1}V_t|+|B_tV_{t_1}|}{R^n}\ll\frac{2^{k+1}}{R}\cdot\frac{(|B_{t_1}|+|B_t|)V_t}{H(A_t,B_t)}.
$$
In other words,
\begin{equation}\label{f_eqx}
v_x\ll
\frac{2^kc^{\frac12}}{R}\cdot\frac{|B_t|+|B_{t_1}|}{\max\{|A_t|^{1/i},|B_t|^{1/j}\}}.
\end{equation}
This means that all pairs $(A_t,B_t)$ under consideration are
situated within some figure defined by \eqref{f_eqx} which we denote
by~$F_x$.  Note that $F_x$ depends on $B_{t_1}$ and $c_1'$ which in
turn is defined by the point $P$, interval $J$ and the class of
intervals $\Delta(L_t)$.   The upshot is that if all lines $L_t$
intersect at one point $P$ and all intervals $\Delta(L_t)$ intersect
$J$ then all pairs $(A_{t},B_{t})$, except possibly one with
$t=t_1$, lie in the set $F_x\cap \LL$.

When considering  inequality~\eqref{dist_y}, similar
arguments enable us to conclude that all pairs $(A_{t},B_{t})$,
except possibly one, lie in the set  $F_y\cap\LL$ where $F_y$ is the
figure defined by
\begin{equation}\label{f_eqy}
v_y\ll
\frac{2^kc^{\frac12}}{R}\cdot\frac{|A_t|+|A_{t_2}|}{\max\{|A_t|^{1/i},|B_t|^{1/j}\}}.
\end{equation}
This together with the previous statement for $F_x$ implies that all
pairs $(A_t,B_t)$, except possibly two, lie in the set $F_x\cap
F_y\cap \LL$.

\subsubsection{Interval $J$
does not satisfy~\eqref{cond_c1d} and intervals $\Delta(L_t)$ are of
Type 1}

Now assume that interval $J$ does not satisfy~\eqref{cond_c1d}. Then
by applying~\eqref{cond_dist2} for the pair of intervals
($\Delta(L_t)$, $\Delta(L_{t_1})$) and ($\Delta(L_t)$,
$\Delta(L_{t_2})$) we obtain the following two conditions:
$$
\frac{|B_{t_1}|\max\{|A_t|,|B_t|\}+|B_t|\max\{|A_{t_1}|,|B_{t_1}|\}}{R^{2n}}\ge
\sigma_x:=\frac{\omega_x(P,J)}{(c'_1)^2 c_x(\vC,i,j)}  \qquad   t\neq
t_1
$$$$
\frac{|A_{t_2}|\max\{|A_t|,|B_t|\}+|A_t|\max\{|A_{t_2}|,|B_{t_2}|\}}{R^{2n}}\ge
\sigma_y:=\frac{\omega_y(P,J)}{(c'_1)^2c_y(\vC,i,j)}  \qquad  t\neq t_2 \, ,
$$
which play the same role as~\eqref{dist_x} and~\eqref{dist_y} in the
previous case. By similar arguments as before, we end
up with two figures $F'_x$ and $F'_y$ defined as follows:
\begin{equation}\label{f_2eqx}
\sigma_x\ll\frac{2^kc^{\frac12}}{R^{n+1}}\cdot\frac{(|B_t|+|B_{t_1}|)\max\{|A_t|,|B_t|,|A_{t_1}|,|B_{t_1}|\}}{V_t\max\{|A_t|^{1/i},|B_t|^{1/j}\}}
\end{equation}
and
\begin{equation}\label{f_2eqy}
\sigma_y\ll\frac{2^kc^{\frac12}}{R^{n+1}}\cdot\frac{(|A_t|+|A_{t_2}|)\max\{|A_t|,|B_t|,|A_{t_2}|,|B_{t_2}|\}}{V_t\max\{|A_t|^{1/i},|B_t|^{1/j}\}}.
\end{equation}
The upshot being that  when $J$ does not
satisfy~\eqref{cond_c1d} all pairs $(A_t,B_t)$, except possibly
two, lie in the set $F'_x\cap F'_y\cap \LL$.

\subsubsection{Intervals $\Delta(L_t)$ are of Type 2}

As usual, for  Type 2 intervals  we assume that~\eqref{cond_tpc1d}
is satisfied.  With appropriate changes, such as the definition of $H(\Delta)$, the same arguments  as above
can be utilised to show that all pairs $(A_t,B_t)$, except possibly two, lie in the
set $F^*_x\cap F^*_y\cap \LL$ where the figures
$F^*_x$ and $F^*_y$ are defined as follows:
\begin{equation}\label{f_tpeqx}
\sigma_x\ll\frac{2^{2k}}{R^2}\cdot\frac{|B_t|}{\max\{|A_t|^{1/i},|B_t|^{1/j}\}}
\end{equation}
and
\begin{equation}\label{f_tpeqy}
\sigma_y\ll\frac{2^{2k}}{R^2}\cdot\frac{|A_t|}{\max\{|A_t|^{1/i},|B_t|^{1/j}\}}.
\end{equation}

\noindent Indeed, the calculations are somewhat simplified since for intervals of Type 2 we have that $|A_t|\asymp
|A_{t_1}|\asymp |A_{t_2}|$ and $|B_t|\asymp |B_{t_1}|\asymp
|B_{t_2}|$.

\subsection{Restrictions to $F_x\cap F_y$ in each class.}

We now use the specific properties of each class to reduce the size of $F_x\cap F_y$ in each case.

\vspace*{2ex}

 $\bullet$ {\bf Class $C(n,k,l,m)$ with $l\ge 1$ and interval $J$ satisfies
\eqref{cond_c1d}.} Consider all intervals $\Delta(L_t(A_t,B_t,C_t))$
from $C(n,k,l,m)$ such that the corresponding coordinates
$(A_t,B_t)$ lie within  the figure  $F_x$ defined by \eqref{f_eqx}. First of all notice that by
\eqref{class_eq1} we have $|A_t|\asymp|B_t|$. Then by
\eqref{class0m_prop} we obtain that
\begin{equation}\label{f_class0_eq1}
\frac{|A_t|}{|V_t|}\asymp 2^m R^{\lambda l}
\end{equation}
which together with~\eqref{asymp_prop}
and~\eqref{f_eqx} implies that
$$
|B_t|\ll \left(\frac{2^k c^{\frac12}}{Rv_x}\right)^{j/i};\quad
|A_t|\ll\left(\frac{2^k c^{\frac12}B_t}{Rv_x}\right)^i\ll
\frac{2^kc^{\frac12}}{Rv_x};\quad V_t\ll
\frac{2^kc^{\frac12}}{Rv_x}2^{-m}R^{-\lambda l}.
$$
If we consider the coordinates $(A_t,B_t)$ within the figure $F_y$
defined by~\eqref{f_eqy},  we obtain the  analogous inequalities:
$$
|A_t|\ll\left(\frac{2^kc^{\frac12}}{Rv_y}\right)^{\frac{i}{j}};\quad
V_t\ll
\left(\frac{2^kc^{\frac12}}{Rv_y}\right)^{\frac{i}{j}}2^{-m}R^{-\lambda
l}.
$$
Hence, it follows that  all coordinates  $(A_t,B_t)\in F_x\cap F_y$
lie inside the box defined by
\begin{equation}\label{f_class0_eq}
|A_t|\ll  \eta :=\min\left\{\frac{2^kc^{\frac12}}{Rv_x},
\left(\frac{2^kc^{\frac12}}{Rv_y}\right)^{\frac{i}{j}}\right\};\quad
|V_t|\ll |A_t|2^{-m}R^{-\lambda l}.
\end{equation}

\vspace*{2ex}

$\bullet$ {\bf Class $C(n,k,l,m)$ with $l\ge 1$ and interval $J$ does not satisfy
\eqref{cond_c1d}.} Consider all intervals $\Delta(L_t(A_t,B_t,C_t))$
from $C(n,k,l,m)$ such that the corresponding coordinates
$(A_t,B_t)$ lie inside~$F'_x$.   As in previous case,  \eqref{f_class0_eq1} is
valid which together with~\eqref{asymp_prop}
and~\eqref{f_2eqx} implies that
$$
|A_t|\ll\left(\frac{2^{k}c^{\frac12}}{\sigma_xR^{n+1}}\cdot
2^mR^{\lambda l}|B_t|\right)^i\ll
\frac{2^{k}c^{\frac12}}{\sigma_xR}\cdot 2^mR^{\lambda l-n};\quad
V_t\ll |A_t|2^{-m}R^{-\lambda l}.
$$
If we consider the coordinates $(A_t,B_t)$ within the figure $F'_y$
defined by~\eqref{f_2eqy},  we obtain the  analogous inequalities:
$$
|A_t|\ll\left(\frac{2^{k}c^{\frac12}}{\sigma_yR}\cdot 2^mR^{\lambda
l-n}\right)^{i/j};\quad V_t\ll |A_t|2^{-m}R^{-\lambda l}.
$$
Denote by $\eta'$ the following minimum
$$
\eta':=\min\left\{\frac{2^{k+m}c^{\frac12}}{\sigma_xR},\left(\frac{2^{k+m}c^{\frac12}}{\sigma_yR}\right)^{i/j}\right\}.
$$
Then,  since for intervals of Type 1 the parameter $l$ is always at
most $l_0$ which in turn satisfies~\eqref{bound_l}, it follows  that all
coordinates $(A_t,B_t)\in F'_x\cap F'_y$ lie inside the box defined
by
\begin{equation}\label{f_class0_2eq}
|A_t|\ll \eta';\quad |V_t|\ll \eta'2^{-m}R^{-\lambda l}.
\end{equation}

\vspace*{2ex}

 $\bullet$ {\bf Class $C^*(n,k)$.} Consider all intervals $\Delta(L_t(A_t,B_t,C_t))$ from $C^*(n,k)$ such that the corresponding coordinates  $(A_t, B_t)\in F^*_x\cap
F^*_y$.   A consequence of that fact that we are considering Type 2 intervals is that  $|B_t|\asymp|B_{t_1}|$. This together with~\eqref{f_tpeqx} and~\eqref{f_tpeqy} implies  that
$$
|B_t|\ll \left(\frac{2^{2k}}{R^2 \sigma_x}\right)^{j/i};\qquad
|A_t|\ll \left(\frac{2^{2k}}{R^2 \sigma_x}|B_t|\right)^i \ll
\frac{2^{2k}}{R^2 \sigma_x};\quad
|V_t|\stackrel{\eqref{type2_eq}}\ll |A_t|R^{-n}
$$
and
$$
|A_t|\ll \left(\frac{2^{2k}}{R^2 \sigma_y}\right)^{i/j}.
$$
Denote by $\eta^*$ the following minimum
$$
\eta^*:=\min\left\{\frac{2^{2k}}{R^2 \sigma_x},\left(\frac{2^{2k}
}{R^2 \sigma_y}\right)^{i/j}\right\}.
$$
The upshot is  that all coordinates $(A_t,B_t)\in F^*_x\cap
F^*_y$ lie inside the box defined by
\begin{equation}\label{f_class0_tpeq}
|A_t|\ll \eta^*;\quad |V_t|\ll \eta^*\cdot R^{-n}.
\end{equation}

\vspace*{2ex}

$\bullet$ {\bf Class  $C_1(n,k)\cap C(n,k,0,m)$.} As
mentioned  in \S\ref{need},  for all subclasses of $C(n,k,0)$,
when consider the intersection with a generic
interval $J$  of length $c_1' R^{-n}$  the constant  $c_1'$   satisfies~\eqref{cond_c1d}.  In other words, $J$ always satisfies~\eqref{cond_c1d}.   With this in mind, consider all intervals $\Delta(L_t(A_t,B_t,C_t))$
from $C_1(n,k)\cap C(n,k,0,m)$ such that the corresponding coordinates
$(A_t,B_t)$ lie within  the figure $F_x$ defined by \eqref{f_eqx}. Then, the  analogue of  \eqref{f_class0_eq1} is
$$
2^m|V_t|\asymp |A_t|  \, .
$$
Although we cannot guarantee that $|B_t|\asymp|B_{t_1}|$,
 by~\eqref{asymp_prop} we have $V_t\asymp V_{t_1}$ and
$|A_t|\asymp|A_{t_1}|$ which in turn implies that
$\max\{|A_t|^{1/i},|B_t|^{1/j}\}\asymp
\max\{|A_{t_1}|^{1/i},|B_{t_1}|^{1/j}\}$. So if $|B_t|\le|B_{t_1}|$, it follows that
$$
\frac{|B_t|+|B_{t_1}|}{\max\{|A_t|^{1/i},|B_t|^{1/j}\}}\asymp
\frac{|B_{t_1}|}{\max\{|A_{t_1}|^{1/i},|B_{t_1}|^{1/j}\}}   \;  .
$$
This together with the previously displayed equation and~\eqref{f_eqx} implies that
$$
|B_t|\le |B_{t_1}|\ll \left(\frac{2^k
c^{\frac12}}{R v_x}\right)^{j/i}.
$$
On the other hand,  if $|B_{t_1}|<|B_t|$ we straightforwardly obtain  the
same estimate for $|B_t|$. So in both cases, we have that
$$|B_t| \ll \left(\frac{2^k
c^{\frac12}}{R v_x}\right)^{j/i}  ; \quad
|A_t|\ll\left(\frac{2^k c^{\frac12}
\max\{|B_t|,|B_{t_1}|\}}{Rv_x}\right)^i\ll
\frac{2^kc^{\frac12}}{Rv_x}  \, ;\quad V_t\ll
\frac{2^kc^{\frac12}}{Rv_x}2^{-m} \, .
$$
If  we consider the coordinates
$(A_t,B_t)$  within  the figure  $F_y$, similar arguments together with inequality~\eqref{f_eqy} yield  the inequalities:
$$
|A_t|\ll\left(\frac{2^kc^{\frac12}}{Rv_y}\right)^{\frac{i}{j}};\quad
V_t\ll
\left(\frac{2^kc^{\frac12}}{Rv_y}\right)^{\frac{i}{j}}2^{-m}R^{-\lambda
l}.
$$
Notice that these inequalities are  exactly the same as when considering `Class $C(n,k,l,m)$ with  $ l\ge 1$, interval $J$ satisfies~\eqref{cond_c1d}' above.  The upshot is that  all coordinates  $(A_t,B_t)\in F_x\cap F_y$
lie inside the box defined by
\begin{equation}\label{f_class1_eq}
|A_t|\ll \eta \, ;\qquad |V_t|\ll 2^{-m} |A_t|   \ .
\end{equation}
Here $\eta$ is as in~\eqref{f_class0_eq} and notice that~\eqref{f_class1_eq}
is indeed equal to~\eqref{f_class0_eq} with $l=0$.

\vspace*{2ex}

$\bullet$ {\bf Class $C_2(n,k)$.} In view of~\eqref{class2_b12}, for  intervals $\Delta(L_t(A_t,B_t,C_t))$ from  $C_2(n,k)$ we have
that $|B_t|\asymp |B_{t_1}|$. Moreover, although we cannot guarantee that $|A_t|\asymp
|A_{t_2}|$,  we still have that  $\max\{|A_t|^{1/i},|B_t|^{1/j}\}\asymp
\max\{|A_{t_1}|^{1/i},|B_{t_1}|^{1/j}\}$  and therefore one can apply the
same arguments as when considering class  $C_1(n,k)\cap C(n,k,0,m)$ above.  As a consequence of~\eqref{f_eqx}
and~\eqref{f_eqy},  it follows that  all coordinates  $(A_t,B_t)\in F_x\cap F_y$
lie inside the box defined by
\begin{equation}\label{f_class2_eq}
|A_t|\ll \eta   \, ; \qquad |B_t|\ll\eta^{j/i}.
\end{equation}

\vspace*{2ex}

$\bullet$ {\bf Class $C_3(n,k,u,v)$.} Consider all intervals $\Delta(L_t(A_t,B_t,C_t))$
from $C_3(n,k,u,v)$ such that the corresponding coordinates
$(A_t,B_t)$ lie within  the figure $F_x$.
In view of~\eqref{class2_b12}, we
have that $|A_t|\asymp |B_t|\asymp |B_{t_1}|\asymp|A_{t_2}|$
and~\eqref{class3_prop} implies that
$\max\{|A_t|^{1/i},|B_t|^{1/j}\}> R^{\lambda u} |B_t|^{1/j}$. This together with~\eqref{f_eqx} implies that
$$
\frac{Rv_x}{2^kc^{\frac12}}\ll B_t^{-i/j}R^{-\lambda u} \Rightarrow
B_t\ll \left(\frac{2^kc^{\frac12}}{Rv_x}\right)^{j/i}R^{-\lambda
uj/i}
$$
and
$$
|A_t|\ll \left(\frac{2^kc^{\frac12}}{Rv_x}|B_t|\right)^i \ll
\frac{2^kc^{\frac12}}{Rv_x}\cdot R^{-\lambda ju}.
$$
If we consider the coordinates $(A_t,B_t)$ within the figure $F_y$
defined by~\eqref{f_eqy},  we obtain the  analogous inequalities:
$$
|A_t|\ll \left(\frac{2^kc^{\frac12}}{Rv_y}\right)^{i/j};\quad
|B_t|\ll \frac{2^kc^{\frac12}}{Rv_y}\cdot R^{-\lambda ju}.
$$
The upshot is that all coordinates $(A_t,B_t)$ from $F_x\cap F_y$ lie inside  the
box defined by
\begin{equation}\label{f_class3_eq}
\begin{array}{l}
\displaystyle|A_t|\ll \,  \eta_3 :=  \min\left\{\frac{2^kc^{\frac12}}{Rv_x}\cdot
R^{-\lambda
ju},\left(\frac{2^kc^{\frac12}}{Rv_y}\right)^{i/j}\right\}\\[5ex]
\displaystyle
|B_t|\ll \,   \eta_3^{j/i}R^{-\lambda ju} =  \min\left\{\left(\frac{2^kc^{\frac12}}{Rv_x}\right)^{j/i}R^{-\frac{\lambda
j}{i}u},\frac{2^kc^{\frac12}}{Rv_y}\cdot R^{-\lambda ju}\right\}.
\end{array}
\end{equation}

\section{The Finale}

The aim of this section is to estimate the number of intervals
$\Delta(L_t)$ from a given  class (either $C(n,k,l,m)$, $C^*(n,k)$,
$C_1(n,k)\cap C(n,k,0,m)$, $C_2(n,k)$ or $C_3(n,k,u,v)$)  that  intersect a fixed  generic interval
$J$ of length $c_1' R^{-n}$. Roughly speaking, the idea is to show that one of the following two situations
necessarily happens:
\begin{itemize}
\item All intervals $\Delta(L_t)$ (except possibly at most two)
intersect the thickening $\Delta(L_0)$ of some line $L_0$.
\item There are not `too many' intervals $\Delta(L_t)$.
\end{itemize}

\noindent As in the previous section we assume that all the
corresponding lines $L_1,L_2,\cdots$ intersect at one point
$P=(p/q,r/q)$.   Then  the quantities $\omega_x(P,J)$ and $\omega_y(P,J)$ are well
defined and the results from \S\ref{sec_inter} are applicable.

\subsection{Point $P$ is  close to $\CCC$}

Assume that
\begin{equation}\label{ineq_omega}
\omega_x(P,J)<\frac{c}{2}q^{-i}\quad\mbox{and}\quad
\omega_y(P,J)<\frac{c}{2}q^{-j}.
\end{equation}
Then, by the definition of $\omega_x$ and $\omega_y$,  we have  that for
each $\Delta(L_t)$
$$
\left|x_t-\frac{p}{q}\right|<\frac{c}{2}q^{-1-i} \quad \mbox{and} \quad
\left|f(x_t)-\frac{r}{q}\right|<\frac{c}{2}q^{-1-j}  \, .
$$
As usual,  $x_t$ is the point in $\Delta(L_t)$ at which $|F'_{L_t}(x)|$ attains its minimum.  In \S\ref{sec_small_dist}, it was shown that  this implies that  all points $x_t$ lie inside $\Delta(L_0)$ for some
line $L_0$. It follows that  all intervals $\Delta(L_t)$ intersect
$\Delta(L_0)$.

\begin{itemize}
\item Assume that $\Delta(L_0)$ has already been removed by the
construction described in \S\ref{tyu}.   In other words,  $\Delta(L_0)\in C(n_0,k_0)$ or
$\Delta(L_0)\in C^*(n_0,k_0)$ with $(n_0,k_0)< (n,k)$). Then each
interval $\Delta(L_t)  \subset \Delta(L_0)$ can be
ignored.   Hence,  the intervals $\Delta(L_t)$ can in total remove at
most two intervals of length
$$
\frac{R}{2^k}\cdot\frac{Kc^{\frac12}}{R^n}
$$
on either side of $\Delta(L_0)$.

\item Otherwise,  by~\eqref{class_eq1} the length of $\Delta(L_0)$ is bounded above
by
$$
\frac{R}{2^k}\cdot\frac{2Kc^{\frac12}}{R^n}.
$$
This implies  that all  the intervals $\Delta(L_t)$ together do not
remove more than a single  interval $\Delta^+(L_0)$ centered at the
same point as $\Delta(L_0)$ but of twice the length. Hence,
the length of the removed interval is bounded above by

\begin{equation}\label{len_rem1}
\frac{R}{2^k}\cdot\frac{4Kc^{\frac12}}{R^n}.
\end{equation}
\end{itemize}

\noindent The upshot is that in either case,  the total length of the intervals removed by
$\Delta(L_t)$ is bounded above  by~\eqref{len_rem1}.

\subsection{Number of intervals $\Delta(L_t)$ intersecting $J$.}\label{subseq2}

We investigate the case when at least one of the bounds in~\eqref{ineq_omega} for
$\omega_x$ or $\omega_y$  is not valid. This implies
the following for the quantities  $v_x$ and $v_y$:
\begin{equation}\label{ineq_vxy}
v_x\ge \frac{cq^{-i}}{2c_1' c_x(\vC,i,j)}\quad\mbox{or}\quad v_y\ge
\frac{cq^{-j}}{2c_1' c_y(\vC,i,j)}.
\end{equation}
The corresponding inequalities  for $\sigma_x$ \and $\sigma_y$ are as
follows:
\begin{equation}\label{ineq_sigxy}
\sigma_x\ge \frac{cq^{-i}}{2(c_1')^2
c_x(\vC,i,j)}\quad\mbox{or}\quad \sigma_y\ge
\frac{cq^{-j}}{2(c_1')^2 c_y(\vC,i,j)}.
\end{equation}
We now  estimate the number of intervals $\Delta(L_t)$ from the same class which
intersect $J$.

A consequence of \S\ref{nicefig}   is that when considering intervals
$\Delta(L_t(A_t,B_t,C_t)$ from the same class  which
intersect $J$, all except possibly at most two of the  corresponding coordinates
$(A_t,B_t)$ lie in the set $F_x\cap F_y\cap \LL$ or $F'_x\cap F'_y\cap \LL$,
or $F^*_x\cap F^*_y\cap \LL$  -- depending on the class  of intervals under consideration. Note that for any two associated lines $L_1 $and $ L_2$,  the coordinates $(A_1,B_1)$,
$(A_2,B_2)$ and $(0,0)$ are not co-linear.    To see this, suppose that the three points did lie on a line. Then
$A_1/B_1 =  A_2/B_2$  and so   $L_1$  and $ L_2$  are parallel.   However, this is impossible since the lines  $L_1$ and $L_{2}$ intersect at the rational point~$P=(p/q,r/q)$.

Now let  $M$  be the number of intervals $\Delta(L_t)$ from the same
class intersecting $J$ and let $F$ denote  the  convex `box' which covers
$F_x\cap F_y$ or $F'_x\cap F'_y$  or $F^*_x\cap F^*_y$ --  depending on the class  of intervals under consideration.
In view of the discussion above, it then follows that  the lattice points of interest in $F\cap \LL$ together with the lattice point $(0,0)$ form
the vertices of  $(M-1)$ disjoint triangles lying within $F$. Since
the area of the fundamental domain of $\LL$ is equal to $q$, the
area of each of these disjoint triangles is at least $q/2$  and
therefore we have that
\begin{equation}\label{m_ineq}
\frac{q}{2}(M-1)  \, \le  \  \area (F).
\end{equation}
We proceed to  estimate $M$  for each class separately.

$\bullet$ {\bf Classes $C(n,k,l,m), l\ge 1$ and $C_1(n,k)\cap C(n,k,0,m)$ and   $J$
satisfies \eqref{cond_c1d}.} By using
either~\eqref{f_class0_eq} for class $C(n,k,l,m), l\ge 1$
or~\eqref{f_class1_eq} for class $C_1(n,k)\cap C(n,k,0,m)$, it follows  that
$$
\area(F) \ \ll \ \eta^2 2^{-m}R^{-\lambda l}
\stackrel{\eqref{ineq_vxy}}\ll
\max\left\{\left(\frac{2^kc_1'}{Rc^{\frac12}}\right)^2,\left(\frac{2^kc_1'}{Rc^{\frac12}}\right)^{2i/j}\right\}
\cdot q^{2i} 2^{-m} R^{-\lambda l}.
$$
This combined with \eqref{m_ineq} gives  the following estimate
\begin{equation}\label{class0_meq}
M   \ \ll \  \max\{D^2, D^{2i/j}\}\cdot2^{-m}R^{-\lambda
l} \ \quad\mbox{where}\quad D:=\frac{2^kc_1'}{Rc^{\frac12}}.
\end{equation}

$\bullet$ {\bf Class $C(n,k,l,m),\; l\ge 1$ and  $J$ does not satisfy
\eqref{cond_c1d}.} By~\eqref{f_class0_2eq},  it follows  that
\begin{eqnarray*}
\area(F)  & \ll  & (\eta')^2 2^{-m}R^{-\lambda l}  \\[2ex]
&\stackrel{\eqref{ineq_sigxy}}\ll  &
\max\left\{\left(\frac{2^{k+m}(c_1')^2}{Rc^{\frac12}}\right)^2,\left(\frac{2^{k+m}(c_1')^2}{Rc^{\frac12}}\right)^{2i/j}\right\}q^{2i}\cdot
2^{-m}R^{-\lambda l}.
\end{eqnarray*}
This combined with \eqref{m_ineq} gives  the following estimate
\begin{equation}\label{class0_2meq}
M  \ \ll \ \max\{(D')^2,(D')^{2i/j}\}   \  2^mR^{-\lambda l} \quad\mbox{where}\quad
D':=\frac{2^k(c_1')^2}{Rc^{\frac12}}.
\end{equation}

$\bullet$  {\bf Class $C^*(n,k)$.} By~\eqref{f_class0_tpeq},  it follows  that
$$
\area(F)  \ \ll \  (\eta^*)^2 R^{-n} \ll
\max\left\{\left(\frac{2^kc_1'}{Rc^{\frac12}}\right)^4,\left(\frac{2^kc_1'}{Rc^{\frac12}}\right)^{4i/j}\right\}
q^{2i}R^{-n} \, .
$$
This combined with \eqref{m_ineq} gives  the following estimate
\begin{equation}\label{class0_tp2meq}
M  \ \ll \  \max\{(D^*)^4, (D^*)^{4i/j}\}\cdot R^{-n}  \quad\mbox{where}\quad
D^*:=\frac{2^kc_1'}{Rc^{1/2}}.
\end{equation}

$\bullet$  {\bf Class $C_2(n,k)$.} By \eqref{f_class2_eq},  it follows that
$$
\area(F) \ \ll  \  \eta^{1+\frac{j}{i}}   \ \ll \ \max\{D^{1/i}, D^{1/j}\}  \, q  \,  .
$$
This combined with \eqref{m_ineq} gives  the following estimate
\begin{equation}\label{class2_meq}
M\ll \max\{D^{1/i},D^{1/j}\}.
\end{equation}

$\bullet$  {\bf Class $C_3(n,k,u,v)$.} By \eqref{f_class3_eq}, it follows that
\begin{eqnarray*}
\area(F)  & \ll &  \eta_3^{1/i} R^{-\lambda uj}  \   \ll  \  \max\left\{(D\cdot
R^{-\lambda uj}q^i)^{1/i}, (D\cdot q^j)^{1/j}\right\}R^{-\lambda uj}
\\[1ex]
& \ll & \max\{ D^{1/i}R^{-\frac{\lambda uj(1+i)}{i}}, D^{1/j}R^{-\lambda
uj}\}\cdot q.
\end{eqnarray*}

\noindent This combined with \eqref{m_ineq} gives  the following estimate
\begin{equation}\label{class3_meq}
M   \ \ll  \  \max\{ D^{1/i}R^{-\frac{\lambda uj(1+i)}{i}},
D^{1/j}R^{-\lambda uj}\}.
\end{equation}

\subsection{Number of subintervals removed by a single interval $\Delta(L)$}

Let $c_1:=c^{\frac12}R^{1+\omega}$  and   $\omega:=ij/4$  be as in  \eqref{ineq_c}.
Consider the nested intervals $J_n\subset J_{n-1}\subset
J_{n-2}\subset\ldots \subset J_0$ where $J_k\in \JJ_k$  with  $0\le k\le
n$.   Consider an interval $\Delta(L)\in
C(n) \cap C^*(n)$ such that $\Delta(L)\cap J_n\neq\emptyset$. We now
estimate the number of intervals
$I_{n+1}\in\II_{n+1}$ such that  $\Delta (L)\cap I_{n+1}
\neq\emptyset$ with $I_{n+1}\subset J_n$.   With reference to the construction of $\JJ_{n+1}$, the desired estimate  is exactly the same as the
number of intervals $I_{n+1}\in\II_{n+1}$ which are removed by the
interval $\Delta(L)$.   By definition, the length of any interval $I_{n+1}$ is
$c_1R^{-n-1}$ and the length of $\Delta(L)$ is
$2Kc^{\frac12}(H(\Delta))^{-1}$. Thus, the number of removed
intervals is bounded above by
\begin{equation}\label{neq4}
2\frac{Kc^{\frac12}R^{n+1}}{c_1H(\Delta)}+2  \ = \ \frac{2KR^{n-\omega}}{H(\Delta)}+2.
\end{equation}
Since  $R^{n-1} \le H(\Delta) < R^n $,  the above  quantity varies between  $2$ and
$[2KR^{1-\omega}]+2$.

\subsection{Condition on $l$ so that  $J_{n-l}$
satisfies \eqref{cond_c1d}}

Consider an interval $J_{n-l}$. Recall, by definition
$$
|J_{n-l}|= c_1R^{-n+l} = (c_1R^{l})\cdot R^{-n}.
$$
So in this case the parameter $c_1'$ associated with the generic interval $J$ is equal to $c_1R^l$ and by
the choice of $c_1$ it clearly   satisfies~\eqref{lbound_c1}.
We now obtain a condition on $l$ so that~\eqref{cond_c1d} is valid when considering  the intersection of intervals from $C(n,k,l,m)$ with $J_{n-l}$.   With this in mind,
on  using the fact that $m \leq \lambda \log_2R$, it follows  that
$$
8\vC \cdot c_1R^{-n+l}\le R^{-\lambda (l+1)} \, \le \,  2^{-m} R^{-\lambda l}  \, .
$$
Thus, \eqref{cond_c1d} is satisfied if
$$
8\vC\cdot
c_1 R^\lambda \cdot R^{l(\lambda+1)}  \, \le  \,  R^n  \, .
$$
By the choice of $c_1$ and in view of~\eqref{ineq_c},  we have  that for $R$ sufficiently large
\begin{equation}\label{ineq_c1}
c_1<\frac{1}{8\vC R^\lambda}.
\end{equation}
Therefore,~\eqref{cond_c1d} is satisfied for $J_{n-l}$  if
$$
l\le\frac{n}{\lambda+1}.
$$
Notice that this is always the case  when $l=0$.

\subsection{Proof of Proposition \ref{prop1}}

Define the parameters $\epsilon:= \frac12(ij)\omega  = \frac18(ij)^2 $ and
\begin{equation}\label{def_ck}
\tilde{c}(k):=\left\{\begin{array}{ll}\displaystyle\frac{c_1R^{\epsilon-\omega}}{2^k}&\mbox{ \ if \  \
}2^k<R^{1-\omega}\\[2ex]
c_1R^{\epsilon-1}&\mbox{ \ if \  \ }2^k\ge R^{1-\omega}.
\end{array}\right.
\end{equation}

\noindent 
Consider an interval $J_{n-l}  \in \JJ_{n-l}$.  Cover $J_{n-l}$ by intervals $J_{l,1},\ldots,
J_{l,d}$ of length $\tilde{c}(k)R^{-n+l}$. Note that by the choice of $c_1$ and $R$ sufficiently large the quantity  $ c_1'=:\tilde{c}(k)R^l $
 satisfies~\eqref{lbound_c1}.   It is easily seen that the number $d$ of such
intervals is estimated as follows:

\begin{equation}\label{d_ineq}
\left\{\begin{array}{ll} d\le 2^kR^{\omega-\epsilon}&\mbox{ \ if \  \
}2^k<R^{1-\omega} \\[2ex]
d\le R^{1-\epsilon}&\mbox{ \ if \  \ }2^k\ge R^{1-\omega}  \, .
\end{array}\right.
\end{equation}


\subsubsection{Part 1 of Proposition \ref{prop1}}

A consequence of  \S\ref{sec_inter} is that if $c_1' = c_1R^l$ satisfies
either \eqref{cond_c1_2} or~\eqref{cond_c1_22}, depending on whether inequality~\eqref{cond_c1d} holds or not , then all lines $L$
associated with intervals $\Delta(L)\in C(n,k,l,m)$ such that
$\Delta(L)\cap J_{n-l}\neq \emptyset$  intersect at a single  point.
This statement remains valid if the interval $J_{n-l}$ is replaced by
any nested interval $J_{l,t}$. Inequality~\eqref{cond_c1_2} is equivalent to
\begin{equation*}
c_1\le \delta\cdot
\left(\frac{2^kc^{\frac12}}{R}R^\lambda\right)^{-\frac{3i}{i+1}}  \quad    \mbox { \ or   \ }  \quad
c^{\frac12} \le \delta\cdot
\left(\frac{2^kc^{\frac12}}{R}R^\lambda\right)^{-\frac{3i}{i+1}}R^{-1-\omega}
\end{equation*}
and inequality~\eqref{cond_c1_22} is equivalent to
$$
c_1\le \delta\cdot
\left(\frac{2^kc^{\frac12}}{R}R^\lambda\right)^{-\frac{i}{i+1}}R^{-\lambda/3}\quad    \mbox { \ or   \ }  \quad
c^{\frac12} \le \delta\cdot
\left(\frac{2^kc^{\frac12}}{R}R^\lambda\right)^{-\frac{i}{i+1}}R^{-1-\omega
-\lambda/3}.
$$
In view of \eqref{ineq_c}, for $R$ large enough both of these  upper bound inequalities on $c$ are satisfied.  Thus
the coordinates $(A,B)$ associated with intervals $\Delta(L(A,B,C))\in
C(n,k,l,m)$ intersecting $J_{l,t}$ where $1\le t\le d$, except possibly at most two, lie within the figure $F:= F_x\cap F_y\cap\LL$ or $F:= F'_x\cap
F'_y\cap\LL$   -- depending on whether or not  $J_{l,t}$ satisfies~\eqref{cond_c1d}. Moreover, note that  the figure $F$ is  the same for  $1\le t\le d$; i.e. it is independent of $t$.

If   \eqref{ineq_omega} is valid,  then all intervals
$\Delta(L)$ that intersect $J_{l,t}$ can remove  at most two
intervals of total length bounded above by
$$
\frac{R}{2^k}\cdot \frac{4Kc^{\frac12}}{R^n}.
$$
Then, it follows that the number of removed intervals
$I_{n+1}\subset J_{n-l}$ is bounded above by
\begin{equation}\label{oneline}
\left(\frac{R}{2^k}\cdot \frac{4Kc^{\frac12}}{R^n} \cdot \frac{1}{|I_{n+1}|}  +4\right)\cdot
d=
4 \left(\frac{K \, R^{1-\omega}}{2^k}+1\right)\cdot
d  \ll  \left(\frac{ \, R^{1-\omega}}{2^k}+1\right)\cdot
d   \stackrel{\eqref{d_ineq}}\ll R^{1-\epsilon}.
\end{equation}

Otherwise, if \eqref{ineq_omega} is false then the number $M$ of
intervals $\Delta(L)\in C(n,k,l,m)$  that intersect some $J_{l,t}$ ($1 \le t \le d$) can be estimated by~\eqref{class0_meq} if  $J_{l,t}$
satisfies~\eqref{cond_c1d} and  by~\eqref{class0_2meq} if~\eqref{cond_c1d} is not satisfied. This leads to the following estimates.

\begin{itemize}
\item {\bf $M$ is bounded by~\eqref{class0_meq} and
$2^k<R^{1-\omega}$.} Then
$$
M\ll \left(\frac{2^k
c^{\frac12}R^{1+\epsilon}R^l}{2^kRc^{\frac12}}\right)^2
2^{-m}R^{-\lambda l}\le (R^\epsilon)^2\cdot
R^{\left(2-\lambda\right)l}.
$$
By~\eqref{cond_lambda}, $\lambda>2$ and therefore $M\ll
R^{2\epsilon}$.
\item {\bf $M$ is bounded by~\eqref{class0_meq} and
$2^k\ge R^{1-\omega}$.} Then
$$
M\ll \left(\frac{2^k
c^{\frac12}R^{\omega+\epsilon}R^l}{Rc^{\frac12}}\right)^2
2^{-m}R^{-\lambda l}\le (R^{\omega+\epsilon})^2.
$$
because $R\ge 2^k$ and $\lambda>2$.
\item {\bf $M$ is bounded by~\eqref{class0_2meq} and
$2^k< R^{1-\omega}$.} Then
$$
M\ll \left(\frac{2^k c
R^{2+2\epsilon}R^{2l}}{2^{2k}Rc^{\frac12}}\right)^2 2^mR^{-\lambda
l}\le  c\cdot \frac{2^m R^2}{2^{2k}} R^{4\epsilon}\cdot
R^{(4-\lambda) l}.
$$
Since $\lambda>4$ by~\eqref{cond_lambda} and $c<R^{-2-\lambda}$
by~\eqref{ineq_c}, it follows that $M\ll R^{4\epsilon}$.
\item {\bf $M$ is bounded by~\eqref{class0_2meq} and
$2^k\ge R^{1-\omega}$.} Then
$$
M\ll \left(\frac{2^k c
R^{2\epsilon+2\omega}R^{2l}}{Rc^{\frac12}}\right)^2 2^mR^{-\lambda
l}\le  c\cdot 2^m \cdot R^{4(\epsilon+\omega)}\cdot R^{(4-\lambda)
l}.
$$
Again,  by the choice of  $\lambda$ and  $c$ it follows that $M \ll R^{4(\epsilon+\omega)}$.
\end{itemize}

\noindent    The upshot of the above upper bounds on $M$ is that
\begin{equation}\label{class1_Meq}
M\ll \left\{\begin{array}{ll}(R^\epsilon)^4&\mbox{if
}2^k<R^{1-\omega}\\[2ex]
(R^{\omega+\epsilon})^4&\mbox{if }2^k\ge R^{1-\omega}.
\end{array}\right.
\end{equation}

\noindent In addition  to these $M$ intervals, we can have at most
another $2d$ intervals -- two for each $1\le t \le d$ corresponding to the fact that there may be up to two  exceptional  intervals $\Delta(L(A,B,C)) $ with associated coordinates  $(A,B)$  lying  outside the figure $F$.  By analogy
with~\eqref{oneline},  these intervals  remove at most $R^{1-\epsilon}$
intervals $I_{n+1}\in\II_{n+1}$ with $I_{n+1}\subset J_{n-l}$.

On multiplying $M$ by the number of intervals
$I_{n+1}\in\II_{n+1}$ removed by each $\Delta(L)$ from $C(n,k,l,m)$,
we obtain  via  \eqref{neq4} that the total number of intervals $I_{n+1}\in \II_{n+1} $ with
$I_{n+1}\subset J_{n-l}$ removed by $\Delta(L)\in C(n,k,l,m)$ is bounded above by
\begin{eqnarray*}
2 \, R^{1-\epsilon}+\left(\frac{2 K R^{n-\omega}}{H(\Delta)}+2\right)\cdot (
R^\epsilon)^4  & \stackrel{\eqref{class_prop1}}\ll & R^{1-\epsilon}+
\left(\frac{R^{1-\omega}}{2^k}+1\right)\cdot(R^\epsilon)^4   \\[2ex]
& \ll  &
R^{1-\epsilon}+R^{1-\omega+4\epsilon}   \quad  {\rm if  }  \quad   2^k<R^{1-\omega}
\end{eqnarray*}
and by
\begin{eqnarray*}
2 \, R^{1-\epsilon}+\left(\frac{2 K R^{n-\omega}}{H(\Delta)}+2\right)\cdot (
R^{\omega+\epsilon})^4 & \stackrel{\eqref{class_prop1}}\ll  &
R^{1-\epsilon}+
\left(\frac{R^{1-\omega}}{2^k}+1\right)\cdot(R^{\epsilon+\omega})^4   \\[2ex]
& \ll  &
R^{1-\epsilon}+R^{4(\omega+\epsilon)}   \quad  {\rm if  }  \quad   2^k \ge R^{1-\omega}  \, .
\end{eqnarray*}
 Since $\omega= \frac14ij$ and
$\epsilon=\frac12 (ij)\omega$, in either case the number of removed
intervals $I_{n+1}$ is  $ \ll~R^{1-\epsilon}$.
Now recall that the  parameters $k$ and $m$ can only take on a constant times $\log R$ values. Hence,  it follows that
$$
\#\{I_{n+1}\in\II_{n+1}: I_{n+1}\subset J_{n-l}, \exists
\Delta(L)\in C(n,l), \Delta(L)\cap I_{n+1}\neq\emptyset\} \ll \log^2
R\cdot R^{1-\epsilon}.
$$
For $R$ large enough the r.h.s.  is bounded above by
$R^{1-\epsilon/2}$.

\subsubsection{Part 2 of Proposition \ref{prop1}}

Consider an interval $J_{n-n_0}  \in \JJ_{n-n_0}$, where  $n_0$ is  defined by~\eqref{def_n0} and $n\ge 3n_0$.   Cover $J_{n-n_0}$ by intervals
$J_{n_0,1},\ldots, J_{n_0,d}$ of length $\tilde{c}(k)R^{-n+n_0}$
where $\tilde{c}(k)$ is defined by~\eqref{def_ck}. Notice that $d$ satisfies~\eqref{d_ineq}. Also,  in view of~\eqref{def_n0}  it follows that  $c_1'
:= \tilde{c}(k)R$  satisfies~\eqref{cond_tpc1d}.
Therefore,  Lemma~\ref{lem_tp2vp} is applicable to  the intervals
$J_{n_0,t}$ with $ 1\le t\le d$ and indeed is applicable to the  whole interval $J_{n-n_0}$.

To ensure that all lines associated with $\Delta(L)\in C^*(n,k)$
such that $\Delta(L)\cap J_{n-n_0}\neq\emptyset$ intersect at one
point, we need to guarantee that~\eqref{cond_c1_tp2} is  satisfied
for $c_1':=c_1R^{n_0}$. This is indeed the case if
\begin{equation}\label{cond_c1_tpm}
c_1R^{n_0}\le \delta\cdot
R^{\frac{j}{1+i}n}\left(\frac{2^k}{R}\right)^{-\frac{2i}{i+1}}.
\end{equation}
Since $i\le j$ we have that  $\frac{j}{1+i}\ge \frac13$ which together with the fact that $n\ge 3n_0$ implies  that  \eqref{cond_c1_tpm} is true if
$$
c^{\frac12} \le \delta\cdot
\left(\frac{2^k}{R}\right)^{-\frac{2i}{i+1}}R^{-1-\omega}
$$
In view of \eqref{ineq_c}, for $R$ large enough  this  upper bound inequality on $c$  is satisfied. Thus the
coordinates $(A,B)$ of all except possibly at most two lines
$L(A,B,C)$ associated with intervals $\Delta(L(A,B,C))\in C^*(n,k)$ with $
\Delta(L)\cap J_{n_0,t}\neq\emptyset$ lie within  the figure
$F:=F^*_x\cap F^*_y\cap \LL$. By analogy with Part 1,
if~\eqref{ineq_omega} is valid then the number of intervals
$I_{n+1}\subset J_{n-n_0} $ removed by intervals  $\Delta(L)$ is bounded above by
$R^{1-\epsilon}$. Otherwise,  the number $M$ of intervals
$\Delta(L)\in C^*(n,k)$ that intersect some $J_{n_0,t}$  ($1\le t \le d $) with associated coordinates  $(A,B)\in F$ can be
estimated by~\eqref{class0_tp2meq}. Thus
$$
M\ll \left(\frac{2^k \tilde{c}(k)}{Rc^{\frac12}}\right)^4 R^{-n}\le
\left\{\begin{array}{lc}
(R^\epsilon)^{4}R^{-n}&\mbox{if }2^k<R^{1-\omega}\\[2ex]
(R^{\epsilon+\omega})^{4}R^{-n}&\mbox{if }2^k\ge R^{1-\omega}.
\end{array}\right.
$$
Since $n\ge 1$ and $\omega+\epsilon<1/4$, it follows that
$M\ll 1$. Now the same arguments as in Part 1 above can be utilized to verify that
$$
\#\{I_{n+1}\in\II_{n+1}: I_{n+1}\subset J_{n-l}, \exists
\Delta(L)\in C^*(n), \Delta(L)\cap I_{n+1}\neq\emptyset\} \ll \log
R\cdot R^{1-\epsilon}  \, .
$$
For $R$ large enough the r.h.s. is bounded above by
$R^{1-\epsilon/2}$.

\subsubsection{Part 3 of Proposition \ref{prop1}}

Consider an interval $J_n  \in \JJ_n$.  Cover $J_n $ by intervals
$J_{0,1},\ldots, J_{0,d}$ of length $\tilde{c}(k)R^{-n}$ where
$\tilde{c}(k)$ is defined by~\eqref{def_ck}.  As before, $d$ satisfies  \eqref{d_ineq}.

First we consider  intervals $\Delta(L)$ from class $C_1(n,k)\cap
C(n,k,0,m)$ such that $\Delta(L)\cap J_n\neq \emptyset$. In this case,  the conditions~\eqref{f_class1_eq} on the convex `box' containing  the
figure $F_x\cap F_y\cap \LL$ and the conditions \eqref{class0_meq}
on $M$ are the same as those when dealing with the  class $C(n,k,l,m)$ in Part~1 above. Thus,  analogous arguments imply that
$$
\#\{I_{n+1}\in\II_{n+1}\;:\; I_{n+1}\subset J_n, \exists
\Delta(L)\in C_1(n), \Delta(L)\cap I_{n+1}\neq\emptyset\}\ll
R^{1-\epsilon/2}.
$$

Next we consider intervals $\Delta(L)$ from  class $C_2(n,k)$  such that $\Delta(L)\cap J_n\neq \emptyset$. A consequence of  \S\ref{sec_inter} is that if $c_1':=c_1$ satisfies
\eqref{cond_c1_3},  then all lines $L$ associated with intervals
$\Delta(L)\in C_2(n,k)$ such that $\Delta(L)\cap J_n\neq \emptyset$
 intersect at a single point. Inequality \eqref{cond_c1_3} is equivalent to
$$
c_1\le\delta\cdot
\left(\frac{2^kc^{\frac12}}{R}R^\lambda\right)^{-1} \quad \mbox{\ or  \ }  \quad
c^{\frac12}\le\delta\cdot
\left(\frac{2^kc^{\frac12}}{R}R^\lambda\right)^{-1}R^{-1-\omega}.
$$
In view of \eqref{ineq_c}, for $R$ large enough  this  upper bound inequality on $c$  is satisfied. Thus the coordinates
$(A,B)$  associated with intervals $\Delta(L(A,B,C))\in C_2(n,k)$
intersecting $J_{0,t} $ where $ 1\le t\le d$, except possibly at most  two, lie within
the figure $F:= F_x\cap F_y\cap\LL$. We now  follow the arguments
from  Part 1.  If~\eqref{ineq_omega} is valid,   then we deduce that  the total
number of intervals $I_{n+1}\subset
J_n$  removed by intervals  $\Delta(L)$ is bounded above by~\eqref{oneline}.  Otherwise,  the number $M$ of intervals
$\Delta(L)\in C_2(n,k)$ that intersect some $J_{0,t}$  ($1\le t \le d $) with associated coordinates  $(A,B)\in F$ can be
estimated by~\eqref{class2_meq}.  Thus, with $c_1' :=
\tilde{c}(k)$ given by~\eqref{def_ck}  we obtain that
$$
M\ll \left(\frac{2^k \tilde{c}(k)}{Rc^{\frac12}}\right)^{1/i} \le
\left\{\begin{array}{lc}
(R^\epsilon)^{1/i}&\mbox{if }2^k<R^{1-\omega}\\[2ex]
(R^{\epsilon+\omega})^{1/i}&\mbox{if }2^k\ge R^{1-\omega}.
\end{array}\right.
$$
It follows  via  \eqref{neq4} that the  total number of intervals $I_{n+1}\in \II_{n+1}$ with $
I_{n+1}\subset J_n$ removed by $\Delta(L)\in C_2(n,k)$ is  bounded above
by
\begin{eqnarray*}
2 \, R^{1-\epsilon}+\left(\frac{2R^{n-\omega}}{H(\Delta)}+2\right)\cdot (
R^\epsilon)^{1/i}   & \stackrel{\eqref{class_prop1}}\ll  &  R^{1-\epsilon}+
\left(\frac{R^{1-\omega}}{2^k}+1\right)\cdot(R^\epsilon)^{1/i}  \\[2ex] & \ll  &
R^{1-\epsilon}+R^{1-\omega+\epsilon/i} \quad  {\rm if  }  \quad   2^k < R^{1-\omega}
\end{eqnarray*}
and
\begin{eqnarray*}
2\, R^{1-\epsilon}+\left(\frac{2R^{n-\omega}}{H(\Delta)}+2\right)\cdot (
R^{\omega+\epsilon})^{1/i} & \stackrel{\eqref{class_prop1}}\ll  &
R^{1-\epsilon}+
\left(\frac{R^{1-\omega}}{2^k}+1\right)\cdot(R^{\epsilon+\omega})^{1/i} \\[2ex] & \ll  &
R^{1-\epsilon}+R^{(\omega+\epsilon)/i} \quad  {\rm if  }  \quad   2^k \ge R^{1-\omega}  \, .
\end{eqnarray*}
 Since $\omega= \frac14ij$ and
$\epsilon=\frac12 (ij)\omega$, in either case the number of removed
intervals $I_{n+1}$ is  $ \ll $ $  R^{1-\epsilon}$.
Hence, we obtain that
$$
\#\{I_{n+1}\in\II_{n+1}\;:\; I_{n+1}\subset J_n, \exists
\Delta(L)\in C_2(n), \Delta(L)\cap I_{n+1}\neq\emptyset\}\ll
 \log
R\cdot R^{1-\epsilon}.
$$
For $R$ large enough the r.h.s.  is bounded above by
$R^{1-\epsilon/2}$.

\subsubsection{Part 4 of Proposition \ref{prop1}}

The proof is pretty much the same as for  Parts 1-3.
Consider an interval $J_{n-u}  \in \JJ_{n-u}$.  Cover $J_{n-u}$ by intervals $J_{u,1},\ldots, J_{u,d}$ of length
$\tilde{c}(k)R^{-n+u}$ where $\tilde{c}(k)$ is given
by~\eqref{def_ck}. As usual,   $d$ satisfies
\eqref{d_ineq}.  Recall, that  $C_3(n,k,u,v)\subset C(n,k,0)$ and for all subclasses of $C(n,k,0)$ when consider the intersection with a generic interval $J$  of length $c_1' R^{-n}$  we require  the constant  $c_1'$  to satisfy~\eqref{cond_c1d}  -- see    \S\ref{need}.   Thus, to begin with we check that  the interval $J_{n-u}$ satisfies~\eqref{cond_c1d}. Now, with $c_1':= c_1 R^u$ and $l=0$, together with the fact that $m \le \lambda \log_2 R$, the desired inequality~\eqref{cond_c1d} would hold if
$$
8\vC c_1R^{u-n}\le R^{-\lambda}  \, .
$$  It
is easily verified that this is indeed true by making use of the
inequalities~\eqref{class3_lbound} and~\eqref{ineq_c1} concerning  $u$ and
$c_1$ respectively.

A consequence of  \S\ref{sec_inter} is that if $c_1':=c_1R^u$ satisfies
\eqref{cond_c1_4},  then all lines $L$ associated with intervals
$\Delta(L)\in C_3(n,k,u,v)$ such that $\Delta(L)\cap J_{n-u}\neq
\emptyset$  intersect at a single  point.  Inequality~\eqref{cond_c1_4} is equivalent to
$$
c_1\le \delta\cdot\frac{R^{1-\lambda
i}}{2^kc^{\frac12}} \quad\mbox{ \ or  \ } \quad c^{\frac12}\le
\delta\cdot\frac{R^{-\lambda i-\omega}}{2^kc^{\frac12}}.
$$
In view of \eqref{ineq_c}, for $R$ large enough this  upper bound inequality on $c$  is satisfied.
Thus  the coordinates
$(A,B)$  associated with intervals $\Delta(L(A,B,C))\in C_3(n,k,u,v)$
intersecting $J_{u,t} $ where $ 1\le t\le d$, except possibly at most  two, lie within
the figure $F:= F_x\cap F_y\cap\LL$.

We now  follow the arguments
from  Part 1.  If~\eqref{ineq_omega} is valid,   then we deduce that  the total
number of intervals $I_{n+1}\subset
J_n$  removed by intervals  $\Delta(L)$ is bounded above by~\eqref{oneline}.  Otherwise,  the number $M$ of intervals
$\Delta(L)\in C_3(n,k,u,v)$ that intersect $J_{u,t}$ with associated coordinates  $(A,B)\in F$ can be
estimated by~\eqref{class3_meq}.  Thus, with $c_1' :=
\tilde{c}(k)$ given by~\eqref{def_ck}  we obtain that
$$
M\ll \left(\frac{2^k
c^{\frac12}R^{1+\epsilon}R^u}{2^kRc^{\frac12}}\right)^{1/i}R^{u(1-\min\{\frac{\lambda
j(1+i)}{i}, \lambda j\})}\stackrel{\eqref{cond_lambda}}\le
(R^\epsilon)^{1/i}   \quad  {\rm if  }  \quad   2^k < R^{1-\omega}
$$
and
$$
M\ll \left(\frac{2^k
c^{\frac12}R^{\omega+\epsilon}R^u}{Rc^{\frac12}}\right)^{1/i}R^{u(1-\min\{\frac{\lambda
j(1+i)}{i}, \lambda j\})}\le (R^{\omega+\epsilon})^{1/i}  \quad  {\rm if  }  \quad   2^k \ge R^{1-\omega}  .
$$

\noindent Note that these are exactly the same estimates for $M$ obtained  in Part 3 above. Then as before, we deduce that the total number of
intervals $I_{n+1}\in \II_{n+1}$ with $  I_{n+1}\subset J_n$ removed by
$\Delta(L)\in C_3(n,k,u,v)$ is bounded above by $R^{1-\epsilon}$. Hence, it follows that
$$
\#\{I_{n+1}\in\II_{n+1}: I_{n+1}\subset J_{n-u}, \exists
\Delta(L)\in \widetilde{C}_3(n,u), \Delta(L)\cap
I_{n+1}\neq\emptyset\} \ll \log^2 R\cdot R^{1-\epsilon}  \, .
$$
For $R$ large enough the r.h.s.  is bounded above by
$R^{1-\epsilon/2}$.

\endproof

\section{Proof of Theorem~\ref{thmnvlines}}
The basic strategy of the proof of Theorem \ref{thmfinite} also
works for Theorem~\ref{thmnvlines}. The key is to establish the
analogue of  Theorem \ref{thmfinitelowerbd}.   In this section we
outline the main differences and  modifications.  Let $(i,j)$ be a
pair of real numbers satisfying \eqref{neq1q}. Given a line
$\L_{\alpha,\beta } : x \to  \alpha x +\beta $ we have that
$$
F_L(x):= (A-B\alpha)x +C-B\beta \qquad {\rm and} \quad  V_L :=
|F'_L(x)| = |A-B\alpha|  \ ,
$$
Thus, with in the context  of Theorem~\ref{thmnvlines} the quantity
$V_L$  is independent of  $x$.  Furthermore,  note that the
Diophantine condition on $\alpha$ implies  that there exists an
$\epsilon > 0  $  such that
\begin{equation}\label{ineq_v_line}
V_L\gg B^{-\frac{1}{i}+\epsilon}.
\end{equation}
Also,  $|F''_L(x)| \equiv 0$  and the analogue of
Lemma~\ref{lem_delta} is the following statement.

\begin{lemma}
There exists an absolute constant $K \geq 1$ dependent  only on
$i,j,\alpha$ and $\beta$ such that
$$
|\Delta(L)|\le K\frac{c}{\maxab\cdot V_L}.
$$
\end{lemma}

A consequence of the lemma is that there are only  Type 1 intervals
to consider.  Next note that for $c$ small enough $H(\Delta)>1$ for
all intervals $\Delta(L)$ . Indeed
$$
H(\Delta) = c^{-1/2}V_L \maxab.
$$
So if $|A|<\frac{|\alpha|}{2} |B|$,  then $V_L\asymp B$ and
$H(\Delta)>1$ follows immediately. Otherwise,
$$
H(\Delta) \stackrel{\eqref{ineq_v_line}}\gg
c^{-1/2}\left(\frac{|A|}{|B|}\right)^{1/i}
$$
which is also greater than 1 for  $c$  sufficiently small.

As in the case of  non-degenerate curves, we partition  the
intervals $\Delta(L)\in \RRR$ into classes $C(n,k,l)$ according to
\eqref{class_prop1} and~\eqref{class_prop2}. Unfortunately,  we can
not guarantee that $\lambda l\le n$ as in the case of curves.
However,  we still have the bound $l\le n$. To see that this is the
case,  suppose  that $l>n$.  Then~\eqref{class_eq1} is satisfied and
\begin{equation}\label{cond_ln}
|V_L|>R^{-\lambda n}(|\alpha|+1)\max\{|A|,|B|\}.
\end{equation}
By~\eqref{class_prop1}, we  have that
$$
R^{n-1}\le H(\Delta)\le R^n.
$$
On combining the previous two displayed  inequalities we get that
$$
|A-\alpha B|\ll |B|^{\frac{i-\lambda}{i(1+\lambda)}}\cdot
(Rc^{1/2})^{\frac{\lambda.}{1+\lambda}}.
$$
Then by choosing $\lambda$ and $c$ such that
\begin{equation}\label{con_lamb3}
\lambda>\frac{i+1}{\epsilon \, i}-1\quad\mbox{ and }\quad
(Rc^{1/2})^{\frac{\lambda.}{1+\lambda}}<\inf_{q\in\NN}\{q^{\frac1i-\epsilon}||q\alpha||\}:=\tau
\end{equation}
implies that
$$
|A-\alpha B|< \tau |B|^{-\frac{1}{i}+\epsilon}.
$$
This contradicts the Diophantine condition imposed on $\alpha$ and
so  we must have that $l \le n$.

\vspace*{3ex}

With the above differences/changes in mind, it is possible to
establish the  analogue of Proposition~\ref{prop1} for lines
$\L_{\alpha,\beta }$ by following the same arguments and ideas as in
the case of $C^{(2)}$ non-degenerate  planar curves.  The key
differences in the analogous  statement for lines is that in Part~1
we have $l\le n$ instead of $\lambda l\le n$ and that Part~2
disappears all together since there are no Type 2 intervals to
consider. Recall, that even when establishing
Proposition~\ref{prop1} for curves, Part 1, 3 and 4 only use the
fact that the curve is two times differentiable -- see  \S\ref{tyu}
Remark 2.   The analogue of Proposition~\ref{prop1} enables us to
construct the appropriate Cantor  set $ \KKK(J_0, \vR,\vr) $ which
in turn leads to  the desired analogue of  Theorem
\ref{thmfinitelowerbd}.

\vspace*{6ex}

\noindent{\bf Acknowledgements.} SV would like to thank Haleh Afshar and  Maurice Dodson   for their fantastic support over the last two decades and for introducing him to Persian rice and Diophantine approximation. What more could a boy possibly want?   SV would also like to thank the one and only Bridget Bennett  for putting up with his baldness and middle age spread.  Finally, much love to those fab girls Ayesha and Iona as they move swiftly into their second decade!


\begin{thebibliography}{99}


\bibitem{Jinpeng}
J.~An.
\newblock Badziahin-Pollington-Velani's theorem and Schmidt's game.
\newblock Preprint  arXiv:1203.2996,  March 2012.


\bibitem{Jinpeng2}
J.~An.
\newblock Two dimensional badly approximable vectors and Schmidt's game
\newblock Preprint  arXiv:1204.3610,  April 2012.




\bibitem{Badziahin-Velani-MAD}
D.~Badziahin,  S.~Velani.
\newblock Multiplicatively badly approximable numbers and generalised {C}antor
  sets.
\newblock {\em Adv. Math.}, 225:2766--2796, 2011.


\bibitem{Badziahin-Pollington-Velani-Schmidt}
D.~Badziahin, A.~Pollington, S.~Velani.
\newblock On a problem in simultaneous diophantine approximation: Schmidt's
  conjecture.
\newblock {\em Ann. of Math. (2)}, 174(3):1837--1883, 2011.


\bibitem{BakerSchmidt-1970}
A.~Baker,  W. Schmidt.
\newblock Diophantine approximation and {Hausdorff} dimension.
\newblock {\em Proc. Lond. Math. Soc.}, 21:1--11, 1970.


\bibitem{Bugeaud-2004}
Y. Bugeaud.
\newblock  {\em Approximation by algebraic numbers. }
\newblock Cambridge Tracts in Mathematics 160, C.U.P., 2004.

\bibitem{Davenport-64:MR0166154}
H.~Davenport.
\newblock A note on {D}iophantine approximation. {II}.
\newblock {\em Mathematika}, 11:50--58, 1964.


\bibitem{DavenportSchmidt-67}
H.~Davenport,  W. Schmidt.
\newblock Approximation to real numbers by quadratic irrationals.
\newblock {\em Acta Arith.}, 13:169--176, 1967.


\bibitem{Kleinbock-Lindenstrauss-Weiss-04:MR2134453}
D.~Kleinbock, E.~Lindenstrauss,  B.~Weiss.
\newblock On fractal measures and {D}iophantine approximation.
\newblock {\em Selecta Math. (N.S.)}, 10(4):479--523, 2004.


\bibitem{Nesharim}
E. Nesharim.
\newblock Badly approximable vectors on a vertical Cantor set.
Preprint   arXiv:1204.0110,   2012.


\bibitem{Pollington-Velani-02:MR1911218}
A.~Pollington,  S.~Velani.
\newblock On simultaneously badly approximable numbers.
\newblock {\em J. London Math. Soc. (2)}, 66(1):29--40, 2002.

\bibitem{Pyartli-1969}
A.~Pyartli.
\newblock Diophantine approximation on submanifolds of euclidean space.
\newblock {\em Funkts. Anal. Prilosz.}, 3:59--62, 1969.
\newblock (In Russian).


\bibitem{Schmidt-1980}
W.~M. Schmidt.
\newblock {\em Diophantine {Approximation}}.
\newblock Springer-Verlag, Berlin and New York, 1980.

\end{thebibliography}
\end{document}